\documentclass[11pt,a4paper]{article}
\usepackage{amsfonts,amsgen,amsmath,amstext,amsbsy,amsopn,amsthm,amsfonts,amssymb,amscd}
\usepackage{epsf,epsfig}
\usepackage{float}
\usepackage{ebezier,eepic}
\usepackage{color}
\usepackage{tikz}
\usepackage{multirow}
\usepackage{graphicx}
\usepackage{subfigure}
\usepackage{mathrsfs}
\usepackage{graphics}
\usepackage{graphicx}
\usepackage{color}
\usepackage{ifpdf}
\usepackage{fancyhdr}
\usepackage{epstopdf}
\usepackage{epsfig}
\usepackage{indentfirst}
\usepackage{CJK}
\newtheorem{thm}{Theorem}[section]

\newtheorem{lem}[thm]{Lemma}
\newtheorem{false statement}{False statement}

\newtheorem{fact}[thm]{Fact}

\theoremstyle{definition}

\newtheorem{claim}{Claim}

\newtheorem{conj}[thm]{Conjecture}

\makeatletter \@addtoreset{equation}{section}

\newcommand{\ex}{{\rm ex}}

\def\hg{\mathcal{G}}

\def\hb{\mathcal{B}}
\begin{document}
\title{A stability result for $C_{2k+1}$-free graphs }
\author{Sijie Ren\footnote{Department of Mathematics, Taiyuan University of Technology, Taiyuan 030024, P. R. China. E-mail:rensijie1@126.com. }\quad\quad
Jian Wang\footnote{Department of Mathematics, Taiyuan University of Technology, Taiyuan 030024, P. R. China. E-mail:wangjian01@tyut.edu.cn. Research supported by Natural Science Foundation of
Shanxi Province No. RD2200004810.} \quad\quad
Shipeng Wang\footnote{Department of Mathematics, Jiangsu University, Zhenjiang, Jiangsu 212013, P. R. China. E-mail:spwang22@yahoo.com. Research supported by NSFC No.12001242.} \quad\quad
Weihua Yang\footnote{Department of Mathematics, Taiyuan University of Technology, Taiyuan 030024, P. R. China. E-mail:yangweihua@tyut.edu.cn. Research supported by NSFC No.11671296. }\\
 }
\date{}
\maketitle

\begin{abstract}
A graph $G$ is called $C_{2k+1}$-free if it does not contain any cycle of length $2k+1$. In 1981, H\H{a}ggkvist, Faudree and Schelp  showed that every $n$-vertex triangle-free graph  with more than $\frac{(n-1)^2}{4}+1$ edges is bipartite.  In this paper,  we extend their result and show that for $1\leq t\leq 2k-2$ and $n\geq 318t^2k$, every $n$-vertex $C_{2k+1}$-free graph with more than $\frac{(n-t-1)^2}{4}+\binom{t+2}{2}$ edges can be made bipartite by either deleting at most $t-1$ vertices or deleting at most $\binom{\lfloor\frac{t+2}{2}\rfloor}{2}+\binom{\lceil\frac{t+2}{2}\rceil}{2}-1$ edges. The construction shows that this is best possible.
\end{abstract}


\section{Introduction}

We consider only simple graphs. For a graph $G$, we use  $V(G)$ and $E(G)$ to  denote its vertex set and  edge set, respectively. Let $e(G)$ denote the number of edges of $G$.  For $S\subseteq{V(G)}$, let $G[S]$ denote the subgraph of $G$ induced by $S$. Let $G-S$ denote the  subgraph induced by $V(G)\backslash{S}$. For simplicity, we write $E(S)$ and $e(S)$ for $E(G[S])$ and $e(G[S])$, respectively. Let $N_G(v)$ denote the set of neighbors of $v$ in $G$. For $S\subseteq V(G)$,  let $N_G(v,S)$ denote the set of neighbors of $v$ in $S$. Let $\deg_G(v)=|N_G(v)|$ and  $\deg(v,S)=|N_G(v,S)|$. The {\it minimum degree} of $G$ is defined as the minimum of $\deg_G(v)$ over all $v\in V(G)$. For a subgraph $H$ of $G$, we also write $N_G(v,H)$, $\deg_G(v,H)$  for $N_G(v,V(H))$, $\deg_G(v,V(H))$, respectively. For $E\subseteq E(G)$, let $G-E$ denote the graph obtained from $G$ by deleting edges in $E$. For disjoint sets $S,T\subseteq V(G)$, let $G[S,T]$ denote the subgraph of $G$ with the vertex set $S\cup T$ and the edge set $\left\{xy\in E(G)\colon x\in S,\ y\in T\right\}$.  Let $e_G(S,T)=e(G[S,T])$.  If the context is clear, we often omit the subscript $G$.

Given a graph $F$, we say that $G$ is {\it $F$-free} if it does not contain $F$ as a subgraph.  The Tur\'{a}n number ${\rm ex}(n,F)$ is defined as  the maximum number of edges in an $F$-free graph on $n$ vertices. Let $T_r(n)$ denote the {\it Tur\'{a}n graph}, the complete $r$-partite graph on $n$ vertices with $r$ partition classes, each of size $\lfloor \frac{n}{r}\rfloor$ or
$\lceil \frac{n}{r}\rceil$. Let $t_r(n)=e(T_r(n))$. The classic Tur\'{a}n theorem~\cite{turan1941} tells us that $T_r(n)$ is the unique graph attaining the maximum number of edges among all $K_{r+1}$-free graphs, where the $r=2$ case is  Mantel's theorem \cite{Mantel}.

An edge $e\in E(H)$ is called a {\it color-critical} edge of $H$ if $\chi(H-e)<\chi(H)$. Let $H$ be a graph with $\chi(H)=r+1$ that contains a  color-critical edge. Simonovits \cite{sim1974} proved that there exists $n_0(H)$ such that if $n>n_0(H)$, then ${\rm ex}(n,H)=t_r(n)$ and $T_r(n)$ is the unique extremal graph. Since $\chi(C_{2k+1})=3$ and $C_{2k+1}$ contains a color-critical edge,  ${\rm ex}(n,C_{2k+1})$ is implied by Simonovits' result for sufficiently large $n$. Dzido \cite{dzido} observed that ${\rm ex}(n,C_{2k+1})$ can be read out from the works of Bondy \cite{bondy1,bondy2}, Woodall \cite{woodall}, and Bollob\'{a}s \cite{bollobas} (pp. 147--156). Later, F\"{u}redi and Gunderson \cite{furedi2015} determined the $ex(n,C_{2k+1})$ for all $n$ and gave a complete description of the extremal graphs.

\begin{thm}[\cite{dzido}, \cite{furedi2015}]\label{thm-1.1}
\label{c2k1bound}
For $k\geq 2$ and $n\geq 4k-2$,
\begin{equation}\label{eq:dirichlet}
\ex(n, C_{2k+1})= \left\lfloor\frac{n^2}{4}\right\rfloor.
\end{equation}
\end{thm}

In 1981, H\H{a}ggkvist, Faudree and Schelp \cite{haggkvist} proved a stability result for  Mantel's Theorem.

\begin{thm}[\cite{haggkvist}]\label{haggkvist}
Let $G$ be a non-bipartite triangle-free graph on $n$ vertices. Then
\[
e(G) \leq \frac{(n-1)^2}{4}+1.
\]
\end{thm}

Let $H_0$ be a graph obtained from $T_2(n-1)$ by replacing an edge $xy$ with a path $xzy$, where $z$ is a new vertex. It is easy to see that $H_0$ is a triangle-free graph with $\frac{(n-1)^2}{4}+1$ edges, which shows that  Theorem \ref{haggkvist} is sharp.

To state our results, we need the following two parameters. Define
\[
d_2(G)=\min \left\{|T|\colon T\subseteq V(G),\ G-T \mbox{ is bipartite}\right\}
\]
and
\[
\gamma_2(G)=\min \left\{|E|\colon E\subseteq E(G),\ G- E \mbox{ is bipartite}\right\}.
\]
Let $H(n,t)$ be a graph obtained from graphs $T_2(n-t-1)$ and $K_{t+2}$ by sharing a vertex. In this paper, we extend Theorem \ref{haggkvist} to the following two results.

\begin{thm}\label{thm-main1}
Let $n,k,t$ be integers with $1\leq t\leq 2k-2$ and $n\geq 318t^2k$. If $G$ is a  $C_{2k+1}$-free ($k\geq 2$) graph on $n$ vertices with $d_2(G)\geq t$, then
$$e(G) \leq \left\lfloor\frac{(n-t-1)^2}{4}\right\rfloor+\binom{t+2}{2},$$
with equality  if and only if $G=H(n,t)$ up to isomorphism.
\end{thm}

\begin{thm}\label{thm-main2}
Let $n,k,t$ be integers with $1\leq t\leq 2k-2$ and $n\geq 252t^2k$. If $G$ is a  $C_{2k+1}$-free ($k\geq 2$) graph on $n$ vertices with $\gamma_2(G)\geq \binom{\lfloor\frac{t+2}{2}\rfloor}{2}+\binom{\lceil\frac{t+2}{2}\rceil}{2}$,  then
\begin{align}\label{ineq-main2}
e(G) \leq \left\lfloor\frac{(n-t-1)^2}{4}\right\rfloor+\binom{t+2}{2},
\end{align}
with equality  if and only if $G=H(n,t)$ up to isomorphism.
\end{thm}

The {\it circumference} $c(G)$ is defined as the length of a longest cycle of $G$. The {\it girth} $g(G)$  is defined as the length of a shortest  cycle of $G$.  We say $G$ is {\it weakly pancyclic} if it contains cycles of all lengths from $g(G)$ to $c(G)$. We need the following result due to Brandt, Faudree and Goddard \cite{brandt}.

\begin{thm}[\cite{brandt}]\label{thm-pancyclic}
Every non-bipartite graph $G$ of order $n$ with minimum degree $\delta(G)\geq\frac{n+2}{3}$ is weakly pancyclic  with girth 3 or 4.
\end{thm}

 Let $P_{\ell+1}$ be a path of length $\ell$ and let $\mathcal{C}_{\geq \ell}$ be the family of cycles with lengths at least $\ell$.  We need the Erd\H{o}s-Gallai Theorem for paths and cycles.

\begin{thm}[\cite{erdos2}]\label{erdos-gallai-1}
For $n\geq \ell\geq 1$,
\[
{\rm ex}(n,P_{\ell+1})\leq \frac{\ell-1}{2}n,\ {\rm ex}(n,\mathcal{C}_{\geq \ell+1})\leq \frac{\ell}{2}(n-1).
\]
\end{thm}
\section{Some useful facts and inequalities}

In this section, we prove some facts and inequalities that are needed.

\begin{fact}\label{fact-1}
Let $C$ be a shortest odd cycle of $G$. If $C$ has length at least 5, then every vertex in $V(G)\setminus V(C)$ has at most two neighbors on $C$.
\end{fact}

\begin{proof}
Let $C=x_0x_1\ldots x_{2m}x_0$ be a shortest odd cycle of $G$, $m\geq 2$. Fix an arbitrary vertex $v\in V(G)\setminus V(C)$.
If $v$ has at least three neighbors on $C$, there exist neighbors $x_i$, $x_j$ of $v$ ($0\leq i<j\leq 2m$) such that the distance between $x_i$ and $x_j$ on $C$ does not equal 2. Then one of $x_ix_{i+1}\ldots x_j vx_i$ and $x_0x_1\ldots x_ivx_jx_{j+1}\ldots x_{2m}x_0$ is an odd cycle shorter than $C$, a contradiction.
\end{proof}

\begin{fact}\label{fact-2}
Let $G$ be a graph with $e(G)> \frac{(n-t-1)^2}{4}$, $t\geq 1$. Let $C$ be a shortest odd cycle of $G$. If $n\geq 6(t+1)$, $C$ has length at most $\frac{2n}{3}$.
\end{fact}

\begin{proof}
 Let $S=V(C)$ and $s=|S|$. Suppose for contradiction that $s\geq \frac{2n}{3}\geq 5$.  By Fact \ref{fact-1}, each vertex in $G-S$ has at most two neighbors on $C$. Then
\begin{align*}
e(G)&=e(S)+e(G-S)+e(S,G-S)\\[5pt]
&\leq s +\binom{n-s}{2}+ 2(n-s)\\[5pt]
&= 2n-s+\frac{(n-s)(n-s-1)}{2}=:f(s).
\end{align*}
Note that $f(s)$  is a decreasing function for $\frac{2n}{3} \leq s\leq n$. Since $n\geq 6(t+1)\geq 12$, it follows that
\[
e(G)\leq f\left(\frac{2n}{3}\right)=\frac{4n}{3}+ \frac{n}{6}\left(\frac{n}{3}-1\right)=\frac{n}{6}\left(\frac{n}{3}+7\right)\leq\frac{n^2}{6}.
\]
Since  $n\geq 6(t+1)$ implies $\frac{n}{12}-\frac{t+1}{2}\geq 0$,
\begin{align*}
e(G) \leq\frac{n^2}{6}\leq \left(\frac{n}{6}+\frac{n}{12}-\frac{t+1}{2}\right)n=\frac{(n-2(t+1))n}{4}<\frac{(n-t-1)^2}{4},
\end{align*}
 a contradiction.
\end{proof}

\begin{fact}\label{fact-3}
Let $G$ be a $C_{2k+1}$-free graph and $C$ be an odd cycle of length $2\ell+1$ in $G$. If $\ell\geq k+1$ then $\deg(x,C) \leq \ell$ for every $x\in V(G)\setminus V(C)$.
\end{fact}

\begin{proof}
Assume that $C=v_0v_1v_2\ldots v_{2\ell}v_0$. Let $x\in V(G)\setminus V(C)$. Define $I_i$ as an indicator function by setting $I_i=1$ for $xv_i\in E(G)$ and setting $I_i=0$ otherwise. Since $G$ is $C_{2k+1}$-free,  one cannot have both $xv_i$ and $xv_{i+2k-1}\in E(G)$ (subscript modulo $2\ell$) for each $i\in \{0,1,\ldots,2\ell\}$.  Hence,
\[
2\deg(x,C)=2\sum_{0\leq i\leq 2\ell} I_i=\sum_{0\leq i\leq 2\ell} (I_i+I_{i+2k-1}) \leq  2\ell+1.
\]
Therefore, $\deg(x,C)\leq \lfloor\frac{2\ell+1}{2}\rfloor =\ell$.
\end{proof}

Let $a, b\geq 0$ be integers. We need the following equalities and inequalities.
\begin{align}
&\left\lceil\frac{a}{2}\right\rceil \left\lfloor\frac{b}{2}\right\rfloor +\left\lfloor\frac{a}{2}\right\rfloor \left\lceil\frac{b}{2}\right\rceil= \left\lfloor\frac{ab}{2}\right\rfloor,\label{equalities2.1-1}\\[5pt]
&\left\lceil\frac{a}{2}\right\rceil \left\lceil\frac{b}{2}\right\rceil +\left\lfloor\frac{a}{2}\right\rfloor \left\lfloor\frac{b}{2}\right\rfloor= \left\lceil\frac{ab}{2}\right\rceil,\label{equalities2.1-2}\\[5pt]
&\left\lfloor\frac{a}{2}\right\rfloor b \leq\left\lfloor\frac{ab}{2}\right\rfloor,\label{equalities2.2}\\[5pt]
&\left\lceil\frac{a}{2}\right\rceil b\leq\left\lceil\frac{ab}{2}\right\rceil+ \left\lfloor\frac{b}{2}\right\rfloor.\label{equalities2.4}
\end{align}
Note that these equalities and inequalities can be easily verified by checking four cases: $a=2i,b=2j$; $a=2i,b=2j+1$; $a=2i+1,b=2j$ and $a=2i+1,b=2j+1$.

Let $a_1,a_2,\ldots,a_m \geq 0$ be integers. By induction on $m$, one can also obtain the following inequality.
\begin{align}
&\left\lceil\frac{a_1+a_2+\cdots+a_m}{2}\right\rceil\geq \left\lceil\frac{a_1}{2}\right\rceil+\left\lfloor\frac{a_2}{2}\right\rfloor
+\cdots+\left\lfloor\frac{a_m}{2}\right\rfloor.\label{equalities2.3}
\end{align}

Define $f(x):= \binom{x}{2}-\left\lfloor\frac{x^2}{4}\right\rfloor+\left\lfloor\frac{x}{2}\right\rfloor$.

\begin{lem}
For integers $a,b,t\geq 0$,
\begin{align}
&f(t+1)=\binom{\lfloor\frac{t+2}{2}\rfloor}{2}+\binom{\lceil\frac{t+2}{2}\rceil}{2},\label{equalities2.6}\\[5pt]
&f(a+b)=f(a)+f(b)+\left\lceil\frac{ab}{2}\right\rceil\label{equalities2.5}.
\end{align}
\end{lem}
\begin{proof}
Note that $\binom{t+1}{2}-\lfloor\frac{(t+1)^2}{4}\rfloor$ can be viewed as the number of edges in $K_{\lfloor (t+1)/2\rfloor}\cup K_{\lceil (t+1)/2\rceil}$. By adding a vertex to $K_{\lfloor (t+1)/2\rfloor}$, we get a graph with $\binom{t+1}{2}-\lfloor\frac{(t+1)^2}{4}\rfloor+ \lfloor (t+1)/2\rfloor=f(t+1)$ edges, which is exactly  a copy of $K_{\lfloor (t+2)/2\rfloor}\cup K_{\lceil (t+2)/2\rceil}$. Thus $\eqref{equalities2.6}$ follows.

Next we prove \eqref{equalities2.5}. Let $A_1, A_2$ be two disjoint vertex sets with  $|A_1|=\left\lceil \frac{a}{2}\right\rceil$, $|A_2|=\left\lfloor \frac{a}{2}\right\rfloor$ and $B_1, B_2$ be two disjoint sets with  $|B_1|=\left\lfloor \frac{b}{2}\right\rfloor$, $|B_2|=\left\lceil \frac{b}{2}\right\rceil$. Clearly $(A_1\cup B_1, A_2\cup B_2)$ is an almost equal partition of $a+b$ vertices. For a set $S$, let $K(S)$ denote the complete graph with vertex set $S$. Note that $\binom{a+b}{2}-\left\lfloor\frac{(a+b)^2}{4}\right\rfloor$ can be viewed as the number of edges in  $K(A_1\cup B_1)\cup K(A_2\cup B_2)$. Similarly, $\binom{a}{2}-\left\lfloor\frac{a^2}{4}\right\rfloor$ can be viewed as the number of edges in $K(A_1)\cup K(A_2)$, $\binom{b}{2}-\left\lfloor\frac{b^2}{4}\right\rfloor$ can be viewed as the number of edges in $K(B_1)\cup K(B_2)$. Thus,
\begin{align*}
f(a+b)-f(a)-f(b)&=|A_1||B_1|+|A_2||B_2|+\left\lfloor \frac{a+b}{2} \right\rfloor -\left\lfloor \frac{a}{2} \right\rfloor-\left\lfloor \frac{b}{2} \right\rfloor\\[5pt]
&=\left\lceil\frac{a}{2} \right\rceil\left\lfloor \frac{b}{2} \right\rfloor+\left\lfloor \frac{a}{2} \right\rfloor\left\lceil \frac{b}{2} \right\rceil+\left\lfloor \frac{a+b}{2} \right\rfloor -\left\lfloor \frac{a}{2} \right\rfloor-\left\lfloor \frac{b}{2} \right\rfloor.
\end{align*}
Define $\delta(a,b):=\lfloor \frac{a+b}{2} \rfloor -\lfloor \frac{a}{2} \rfloor-\lfloor \frac{b}{2} \rfloor$.
Note that $\delta(a,b)$ equals 1 if $a,b$ are both odd and equals 0 otherwise. By \eqref{equalities2.1-1},
\[
f(a+b)-f(a)-f(b)=\left\lfloor\frac{ab}{2}\right\rfloor+\delta(a,b) = \left\lceil\frac{ab}{2}\right\rceil.
\]
\end{proof}

\section{Two structure lemmas}

\begin{lem}\label{lem-2.1}
Let $G$ be a $C_{2k+1}$-free ($k\geq 2$) graph on $n$ vertices with $e(G)\geq \lfloor\frac{(n-t-1)^2}{4}\rfloor+\binom{t+2}{2}$, $t\geq 1$. If $n\geq 78t^2k$, then there exists $T\subseteq V(G)$ with $|T| \leq 3(t+1)$ such that $G-T$ is a bipartite graph on partite sets $X,Y$ and (i), (ii), (iii), (iv) hold.
\begin{itemize}
  \item[(i)] $\delta(G-T)\geq \frac{n-3t-1}{3}$;
  \item[(ii)] $\frac{n}{2}-2\sqrt{tn}\leq |X|,|Y|\leq \frac{n}{2}+\sqrt{3tn}$;
  \item[(iii)] $|\left\{u\in V(G)\setminus T\colon \deg_{G-T}(u)\leq \frac{n}{2}-\sqrt{3tn}\right\}|\leq \sqrt{tn}$;
  \item[(iv)] For $v\in T$,  either $N(v,X\cup Y)\subset X$ or $N(v,X\cup Y)\subset Y$ holds.
\end{itemize}
\end{lem}
\begin{proof}
Set $G_0:=G$.  We obtain a bipartite subgraph from $G_0$ by a vertex-deletion  procedure. Let us start with $i=0$ and obtain a sequence of graphs $G_1,G_2,\ldots,G_\ell,\ldots$ as follows: In the $i$th step, if there exists a  vertex $v_{i+1}\in V(G_i)$ with $\deg_{G_i}(v_{i+1})<\frac{(n-i)+2}{3}$ then let $G_{i+1}$ be the graph obtained from $G_i$ by deleting $v_{i+1}$ and go to the $(i+1)$th step. Otherwise we stop and set $G^*=G_i$,  $T=\{v_1,v_2,\ldots,v_i\}$.  Clearly, the procedure will end  in finite steps.

%First we show that $G^*$ is not empty. Indeed, otherwise for $n\geq 4(t+3)$,
%\begin{align*}
%e(G)< \sum_{i=1}^n \frac{(n-i)+2}{3} =\frac{n(n+3)}{6}
%&=\frac{(n-t-1)^2}{4}+\frac{n}{2}\left(t+2-\frac{n}{6}\right)
%-\frac{(t+1)^2}{4}\\[5pt]
%&\leq \frac{(n-t-1)^2}{4}+2(t+3)\left(t+2-\frac{2(t+3)}{3}\right)
%-\frac{(t+1)^2}{4}\\[5pt]
%&< \frac{(n-t-1)^2}{4}+\binom{t+2}{2}-1,
%\end{align*}
%a contradiction.
%Thus $G^*$ is not empty.

Let $n'=|V(G^*)|$. We claim that $n'\geq 4k$. Indeed, otherwise
\begin{align*}
e(G)< \sum_{i=0}^{n-4k-1}\frac{(n-i)+2}{3}+\binom{4k}{2}&=\frac{(n+4k+5)(n-4k)}{6}+2k(4k-1)\\[5pt]
&< \frac{(n+7k)(n-4k)}{6}+8k^2\\[5pt]
&<\frac{(n-t-1)^2}{4}+\frac{(t+1)n}{2}-\frac{n^2}{12}+\frac{kn}{2}+4k^2.
\end{align*}
By $t+1<2k$, we have $e(G)<\frac{(n-t-1)^2}{4}-\frac{n^2}{12}+\frac{3}{2}kn+4k^2$.  Since $n\geq 78t^2k> 36k$ implies that $\frac{3}{2}kn<\frac{n^2}{24}$, we have $\frac{n^2}{12}>\frac{3}{2}kn+4k^2$. Thus $e(G)<\frac{(n-t-1)^2}{4}$, a  contradiction. Thus $n'\geq 4k$. Suppose to the contrary that $|T|\geq 3t+4$.  Since $G^*$ is $C_{2k+1}$-free, by  Theorem \ref{thm-1.1} we have $e(G_{3t+4})\leq \frac{(n-3t-4)^2}{4}$. Thus,
\begin{align*}
e(G)&< \sum_{i=0}^{3t+3} \frac{(n-i)+2}{3} + \frac{(n-3t-4)^2}{4}\\[5pt]
&\leq\frac{(3t+4)n}{3}+ \frac{(n-3t-4)^2}{4}\\[5pt]
&=\frac{(3t+4)n}{3}+\frac{(n-t-1)^2}{4}+\frac{(t+1)n-(3t+4)n}{2}+\frac{(3t+4)^2-(t+1)^2}{4}\\[5pt]
&=\frac{(n-t-1)^2}{4}-\frac{n}{6}+\frac{(2t+3)(4t+5)}{4}\\[5pt]
&<\frac{(n-t-1)^2}{4}-\frac{n}{6}+2(t+2)^2.
\end{align*}
As $n\geq 78t^2k \geq 12(t+2)^2$, we get $e(G)<\frac{(n-t-1)^2}{4}$,
 a contradiction. Thus $|T|\leq 3t+3$ and  (i) follows from $\delta(G^*)\geq \frac{n'+2}{3}$.

 Suppose that $G^*$ is non-bipartite. As $\delta(G^*)\geq \frac{n'+2}{3}$, by Theorem \ref{thm-pancyclic} we infer that $G^*$ is weakly pancyclic with girth 3 or 4. By (i) and $n\geq 78t^2k>6k+3t+1$,  $G^*$ contains a cycle of length $\delta(G^*)+1\geq \frac{n-3t+2}{3}\geq 2k+1$. Then  the weakly pancyclicity implies that  $G^*$ contains a cycle of length $2k+1$, a contradiction. Thus $G^*$ is bipartite.

Let $X$, $Y$ be two partite sets of $G^*$. We are left to show (ii), (iii) and (iv).  Note that
\begin{align*}
|X||Y|\geq e(G^*)&\geq e(G)-\sum_{i=0}^{3t+2} \frac{(n-i)+2}{3}\\[5pt]
&\geq
\frac{(n-t-1)^2}{4}-\frac{(3t+3)n}{3}\\[5pt]
&\geq \frac{n^2}{4}-\frac{3}{2}(t+1)n\\[5pt]
&\geq \frac{n^2}{4}-3tn.
\end{align*}
Note that $|X|+|Y|=n-|T|$. Set $|X| =\frac{n-|T|}{2}+d$ and $|Y| =\frac{n-|T|}{2}-d$. Then
\[
\frac{(n-|T|)^2}{4}-d^2 = |X||Y| \geq \frac{n^2}{4}-3tn\geq \frac{(n-|T|)^2}{4}-3tn.
\]
It follows that $d\leq \sqrt{3tn}$. Thus,
\[
\frac{n-|T|}{2}-\sqrt{3tn}  \leq |X|,|Y|\leq \frac{n-|T|}{2}+\sqrt{3tn} \leq \frac{n}{2}+\sqrt{3tn}.
\]
By $n\geq 78t^2k> \frac{9}{(2-\sqrt{3})^2}t$ we have $(2-\sqrt{3})\sqrt{tn}\geq 3t$. Then
\[
|X|,|Y|\geq  \frac{n-3t-3}{2}-\sqrt{3tn}\geq\frac{n}{2}-3t-\sqrt{3tn}\geq \frac{n}{2}-2\sqrt{tn}.
\]
Thus (ii) holds.

To show (iii), suppose indirectly that there are at least $\sqrt{tn}$ vertices in $G^*$ with  degree less than $\frac{n}{2}-\sqrt{3tn}$. Let $U\subseteq V(G)$ with $|U|=\lfloor\sqrt{tn}\rfloor$ such that $\deg_{G^*}(u)\leq  \frac{n}{2}-\sqrt{3tn}$ for all $u\in U$. Then  $\deg_G(u)\leq \frac{n}{2}-\sqrt{3tn}+|T|$ for all $u\in U$.
Note that $n\geq 78t^2k>4t$ implies $\sqrt{tn}\leq \frac{n}{2}$ and thereby $n-|U|\geq n-\sqrt{tn} \geq \frac{n}{2}\geq 4k$.
Since $G-U$ is $C_{2k+1}$-free, using Theorem \ref{thm-1.1} we get $e(G-U)\leq \frac{(n-|U|)^2}{4}$. Thus, by $|T|\leq 3(t+1)$ we obtain that
\begin{align*}
e(G)&\leq e(G-U)+|U|\left( \frac{n}{2}-\sqrt{3tn}+|T|\right)\\[5pt]
&= \frac{n^2}{4}+\frac{|U|^2}{4}-\sqrt{3tn}|U|+|T||U|\\[5pt]
&\leq \frac{n^2}{4}+\frac{|U|^2}{4}-\sqrt{3}|U|^2+3(t+1)\sqrt{tn}\\[5pt]
&\leq \frac{n^2}{4}-\left(\sqrt{3}-\frac{1}{4}\right)tn+3(t+1)\sqrt{tn}\\[5pt]
&\leq \frac{(n-t-1)^2}{4}-\left(\sqrt{3}-\frac{5}{4}\right)tn+3(t+1)\sqrt{tn}.
\end{align*}
By $n\geq 78t^2k\geq \frac{36}{(\sqrt{3}-\frac{5}{4})^2}t$ we have $\left(\sqrt{3}-\frac{5}{4}\right)tn\geq 6t\sqrt{tn}\geq3(t+1)\sqrt{tn}$. Hence $e(G)<\frac{(n-t-1)^2}{4}$, a contradiction. Therefore (iii) holds.

To show (iv), suppose for contradiction that for some $v\in T$ there exist $x\in N(v,X)$ and $y\in N(v, Y)$.  Let $X_1=N(y,X)\setminus \{x\}$ and $Y_1=N(x,Y)\setminus \{y\}$.  Note that
\[
|X_1|,|Y_1|\geq \delta(G-T) -1\geq \frac{n-3t-1}{3}-1>  \frac{n}{3}-3t.
\]
Since $G$ is $C_{2k+1}$-free and $k\geq 2$, $G[X_1\cup Y_1]$ contains no  paths of length $2k-3$. By Theorem \ref{erdos-gallai-1} we infer $e(X_1,Y_1) \leq \frac{(2k-3)-1}{2}(|X_1|+|Y_1|)\leq (k-2)n$. Thus,
\begin{align*}
e(G)&\leq e(G-T) +\sum_{0\leq i< |T|} \frac{n-i+2}{3} \\[5pt]
&\leq |X||Y|-|X_1||Y_1|+e(X_1,Y_1)+\frac{n+2}{3}|T|\\[5pt]
&\leq \frac{(n-|T|)^2}{4}-\left( \frac{n}{3}-3t\right)^2+(k-2)n+\frac{n+2}{3}|T|.
\end{align*}
Let $f(x):=\frac{(n-x)^2}{4}+ \frac{n+2}{3}x$. It is easy to check that $f(x)$ is a decreasing function on $[1,\frac{n}{3}-2]$.
Since $n\geq 78t^2k>9t+15$ implies $|T| \leq 3(t+1)\leq \frac{n}{3}-2$,  $f(|T|)\leq f(1)=\frac{(n-1)^2}{4}+\frac{n+2}{3}$. Hence,
\begin{align*}
e(G)&\leq \frac{(n-1)^2}{4} -\left( \frac{n}{3}-3t\right)^2+(k-2)n+\frac{n+2}{3}\\[5pt]
&\leq\frac{(n-t-1)^2}{4}+\frac{tn}{2}-\frac{n^2}{9}+2tn+kn\\[5pt]
&\leq \frac{(n-t-1)^2}{4}- \left(\frac{n^2}{9}- \left(k+\frac{5}{2}t\right)n\right).
\end{align*}
By $n\geq 78t^2k>9k+13t$ we see that $\frac{n^2}{9}\geq (k+\frac{5}{2}t)n$. Then $e(G)\leq \frac{(n-t-1)^2}{4}$,
a contradiction. Thus (iv) holds and the lemma is proven.
\end{proof}

A vertex $u$ is called a {\it high-degree} vertex if $\deg(u)\geq 2\sqrt{tn}$,  otherwise it is called a {\it low-degree} vertex. Let $B(G)$ be the set of all low-degree vertices in $G$, that is,
\[
B(G)=\left\{u\in V(G)\colon \deg(u)< 2\sqrt{tn}\right\}.
\]

\begin{lem}\label{lem-2.2}
Let $G$ be a $C_{2k+1}$-free ($k\geq 2$) graph on $n$ vertices  with $e(G)\geq \lfloor\frac{(n-t-1)^2}{4}\rfloor+\binom{t+2}{2}$, $t\geq 1$ and let $B=B(G)$. If $n\geq 252t^2k$, then
$G-B$ is a bipartite graph on bipartite sets $X,Y$ with (i)-(vi) hold.
\begin{itemize}
  \item[(i)] $|B|\leq t+1$.
  \item[(ii)] $\frac{n}{2}-2\sqrt{tn}\leq |X|,|Y|\leq \frac{n}{2}+2\sqrt{tn}$.
  \item[(iii)] $|\left\{v\in V(G)\setminus B\colon \deg_{G-B}(v)\leq \frac{n}{2}-\sqrt{3tn}\right\}|\leq \frac{3\sqrt{tn}}{2}$.
  \item[(iv)] For $v\in B$,  either $N(v,X\cup Y)\subset X$ or $N(v,X\cup Y)\subset Y$ holds.
  \item[(v)] For $uv\in E(B)$, there does not exist $x_1,x_2\in X$ such that $ux_1,vx_2\in E(G)$. Similarly,  there does not exist $y_1,y_2\in Y$ such that $uy_1,vy_2\in E(G)$.
      \item[(vi)] Let $uwv$ be a path in $B$. If $\deg(w, X\cup Y)\leq 1$, then there does not exist $x\in X, y\in Y$ such that $uy,vx\in E(G)$ or $ux,vy\in E(G)$.
\end{itemize}
\end{lem}
\begin{proof}
By Lemma \ref{lem-2.1} there exists $T\subseteq V(G)$ with $|T| \leq 3(t+1)$ such that $G-T$ is a bipartite graph and (i)-(iv) of Lemma \ref{lem-2.1} hold. Let $X_0$, $Y_0$ be two partite sets of $G-T$. Since $n\geq 252t^2k>144t$, by Lemma \ref{lem-2.1} (i) we have $\delta(G-T)\geq \frac{n-t-1}{3}> \frac{n}{6}\geq  2\sqrt{tn}$. Thus $B\subset T$.

Suppose indirectly that  $G-B$ is not bipartite. Let $C=v_0v_1v_2\ldots v_{2\ell} v_0$ be a shortest odd cycle in $G-B$  and let $S=V(C)$. Clearly, $v_0,v_1,\ldots, v_{2\ell}$ are  all high-degree vertices. Since $C$ is a shortest odd cycle,  $\deg(v_i,C)= 2$ for each $v_i\in S$, that is, $C$ has no chord.  Since $n\geq 252t^2k>\frac{(3t+7)^2}{t}$ implies  $\sqrt{tn}\geq 3t+7\geq |T|+4$,
\[
|N(v_i, (X_0\cup Y_0)\setminus S)|=\deg(v_i)-|T|-\deg(v_i,C)\geq 2\sqrt{tn}-|T|-2\geq \sqrt{tn}+2.
\]
By Lemma \ref{lem-2.1} (iii), each $v_i$ in $S$ has at least two neighbors in $X_0\cup Y_0\setminus S$ with degree at least $\frac{n}{2}-\sqrt{3tn}$.
%Thus $|N^*(v_i,X\cup Y)|\geq 2$ for $i=1,2,\ldots,\ell$.
By Lemma \ref{lem-2.1} (iv), for each $v_i$ in $S$ either $N(v_i,X_0\cup Y_0\setminus S)\subset X_0\setminus S$ or $N(v_i,X_0\cup Y_0\setminus S)\subset Y_0\setminus S$ holds. As $C$ is an odd cycle,  there exist $v_i$, $v_{i+1}$ on $C$ such that $N(v_i,X_0\cup Y_0\setminus S)$  and $N(v_{i+1},X_0\cup Y_0\setminus S)$ are contained in the same one of $X_0\setminus S$ and $Y_0\setminus S$. Without loss of generality, let $x_1, x_2$ be two distinct vertices in $X_0\setminus S$  such that  $v_ix_1,v_{i+1}x_2\in E(G)$ and  $\deg_{G-T}(x_1), \deg_{G-T}(x_2)\geq \frac{n}{2}-\sqrt{3tn}$. Let $Y_1=N(x_1,Y_0)\cap N(x_2,Y_0)\setminus \{v_i, v_{i+1}\}$. By Lemma \ref{lem-2.1} (ii) we know  $|Y_0|\leq \frac{n}{2}+\sqrt{3tn}$. Then
\begin{align*}
|Y_1|&\geq  |N(x_1,Y_0)|+|N(x_2,Y_0)|-|N(x_1,Y_0)\cup N(x_2,Y_0)|-2\\
&\geq \deg_{G-T}(x_1)+\deg_{G-T}(x_2)-|Y_0|-2\\
&\geq\frac{n}{2}-3\sqrt{3tn}-2.
\end{align*}
Let $X_1=X_0\setminus \{x_1,x_2,v_i, v_{i+1}\}$. Note that
\begin{align*}
\delta(G[X_1\cup Y_1]) \geq \delta(G-T)-(|Y_0|-|Y_1|)&\geq \frac{n-3t-1}{3}-4\sqrt{3tn}-2\\[5pt]
&\geq \frac{n}{3}-t-4\sqrt{3tn}-3.
\end{align*}
By $n\geq 252t^2k\geq 504t$ we infer that
\[
\frac{n}{3}-t-4\sqrt{3tn}\geq \frac{n}{3}-\frac{n}{504}-4\sqrt{3tn}\geq \frac{167n}{504}-4\sqrt{3tn}.
\]
As $n\geq 504t$ implies $4\sqrt{3tn}\leq\frac{4n}{\sqrt{168}}< \frac{163n}{504}$, we have $\frac{167n}{504}-4\sqrt{3tn} \geq \frac{n}{126}\geq 2k$. Hence $\delta(G[X_1\cup Y_1])\geq 2k-3$. Then there is a path of length at least $\delta(G[X_1\cup Y_1])+1\geq 2k-2$. For $k\geq 3$ we  choose a  path $P$ of length $2k-4$ with both ends in  $Y_1$. For $k=2$ we simply choose $P$ as  a vertex in $Y_1$. Then $P$ together with $x_1v_iv_{i+1}x_2$ form a cycle of length $2k+1$, a contradiction. Thus $G-B$ is bipartite.

Now we consider (i) and suppose indirectly that $|B|\geq t+2$.
Since $G-B$ is bipartite,  $e(G-B)\leq \frac{(n-|B|)^2}{4}$. It follows that
\[
e(G)< e(G-B)+2\sqrt{tn}|B|\leq \frac{(n-|B|)^2}{4} +2\sqrt{tn} |B|.
\]
Let $f(x)=\frac{(n-x)^2}{4} +2\sqrt{tn} x$.
Since $f(x)$ is convex and $t+2\leq |B|\leq n$,
\[
e(G)=f(|B|)\leq \max\{f(t+2),f(n)\}.
\]
Note that
\begin{align*}
f(t+2)=& \frac{(n-t-2)^2}{4}+2(t+2)\sqrt{tn} \\[5pt]
&\leq \frac{(n-t-1)^2}{4} -\frac{n-t-2}{2}+2(t+2)\sqrt{tn} \\[5pt]
&<\frac{(n-t-1)^2}{4}+\binom{t+2}{2}-1-\frac{n}{2}+2(t+2)\sqrt{tn}.
\end{align*}
 Using $n\geq 252t^2k> 16(t+2)^2t$, we have $\frac{n}{2}\geq 2(t+2)\sqrt{tn}$. Thus $f(t+2)<\frac{(n-t-1)^2}{4}+\binom{t+2}{2}-1$. By $n\geq 252t^2k>144t$ we also have $2\sqrt{tn} \leq \frac{n}{6} \leq  \frac{n-4t}{4} \leq \frac{n-2(t+1)}{4}$. It implies that
\begin{align*}
f(n)=2n\sqrt{tn} \leq \frac{(n-2(t+1))n}{4} <  \frac{(n-t-1)^2}{4}.
\end{align*}
Thus  $e(G)< \frac{(n-t-1)^2}{4}+\binom{t+2}{2}-1$, a contradiction. Therefore (i) holds.

Let $X$, $Y$ be two partite sets of $G-B$. By $B\subseteq T$ we infer $X_0\subset X$ and $Y_0\subset Y$.  By Lemma \ref{lem-2.1} (ii),  $|X|\geq |X_0|\geq \frac{n}{2}-2\sqrt{tn}$ and $|Y|\geq |Y_0|\geq \frac{n}{2}-2\sqrt{tn}$. By $|X|\leq \frac{n}{2}+\sqrt{3tn}$ and $|T\setminus B| \leq |T|\leq 3t+3$ we have
\[
|X|\leq |X_0|+|T\setminus B|\leq  \frac{n}{2}+\sqrt{3tn}+3t+3\leq \frac{n}{2}+\sqrt{3tn}+6t.
\]
By  $n\geq 252t^2k\geq \frac{36}{(2-\sqrt{3})^2}t$, we have $(2-\sqrt{3})\sqrt{tn}\geq 6t$. It
 implies that $|X|\leq \frac{n}{2}+2\sqrt{tn}$.
Similarly, $|Y|\leq \frac{n}{2}+2\sqrt{tn}$. Thus (ii) holds.

Since $B\subseteq T$, $G-T$ is a subgraph of $G-B$. For $n\geq 144t$ we have $\frac{\sqrt{tn}}{2}\geq 6t\geq 3t+3$. By Lemma \ref{lem-2.1} (iii) we infer that
\begin{align*}
\left|\left\{v\in V(G)\setminus B\colon \deg_{G-B}(v)\leq \frac{n}{2}-\sqrt{3tn}\right\}\right| \leq \sqrt{tn}+|T\setminus B| \leq\sqrt{tn}+3t+3 \leq \frac{3\sqrt{tn}}{2}.
\end{align*}
Thus (iii) holds.

We are left to show (iv), (v) and (vi).  Suppose indirectly that  there exists $v\in B$ such that $v$ has two neighbors $x\in X$ and $y\in Y$. By the definition of $B$, we have $\deg(x),\deg(y) \geq 2\sqrt{tn}$. Let $X_1=N(y,X)\setminus \{x\}$ and $Y_1=N(x,Y)\setminus \{y\}$. Then
\[
|X_1|\geq \deg(x)-|B|-1\geq 2\sqrt{tn}-t-2\geq 2\sqrt{tn}-3t.
\]
Since $n\geq 252t^2k> 64t$ implies that $\frac{\sqrt{tn}}{2}\geq 4t \geq 3t+1$, we have $ |X_1|\geq\frac{3\sqrt{tn}}{2}+1$. Similarly $|Y_1|\geq\frac{3\sqrt{tn}}{2}+1$.

If $k=2$, then the $C_5$-free property implies that there is no edge between $X_1$ and $Y_1$. It follows that
$e(G-B)\leq |X||Y|-|X_1||Y_1|\leq \frac{(n-|B|)^2}{4} -tn$.
Then
\begin{align*}
e(G)&\leq e(G-B) +2\sqrt{tn}|B|\leq \frac{(n-|B|)^2}{4}-tn+2\sqrt{tn}|B|.
\end{align*}
Let $f(x):=\frac{(n-x)^2}{4}+2\sqrt{tn}x$. It is easy to check that $f(x)$ is decreasing on $[1,n-4\sqrt{tn}]$. Since $n\geq 252t^2k>25t$ implies $4\sqrt{tn}\leq \frac{4n}{5}$, we have $ |B| \leq t+1\leq n-4\sqrt{tn}$. Then
\begin{align*}
e(G)\leq f(|B|)-tn\leq f(1)-tn=\frac{(n-1)^2}{4}+2\sqrt{tn}-tn\leq \frac{(n-t-1)^2}{4}-\frac{tn}{2}+2\sqrt{tn}.
\end{align*}
By $n\geq 252t^2k$ we have $\frac{tn}{2}\geq 2\sqrt{tn}$. It follows that  $e(G)\leq\frac{(n-t-1)^2}{4}$,
a contradiction. Thus  we may assume that $k\geq 3$.

Since $|X_1|$, $|Y_1|\geq \frac{3\sqrt{tn}}{2}+1$, by (iii) there are $x_1\in X_1$ and $y_1\in  Y_1$ such that $\deg(x_1,Y), \deg(y_1,X)\geq\frac{n}{2}-\sqrt{3tn}$. Let $X_2=N(y_1, X)\setminus \{x,x_1\}$, $Y_2=N(x_1, Y)\setminus \{y,y_1\}$. Then
\[
|X_2|\geq \deg(y_1,X)-2\geq \frac{n}{2}-\sqrt{3tn}-2,\quad |Y_2|\geq \deg(x_1,Y)-2\geq \frac{n}{2}-\sqrt{3tn}-2.
\]
Since $G$ is $C_{2k+1}$-free and $k\geq 3$, $G[X_2\cup Y_2]$ contains no  paths of length $2k-5$. By Theorem \ref{erdos-gallai-1} we infer $e(X_2,Y_2) \leq \frac{(2k-5)-1}{2}(|X_2|+|Y_2|)\leq (k-3)n$. Thus,
\begin{align*}
e(G)&\leq e(G-B) +2\sqrt{tn}|B| \\[5pt]
&\leq |X||Y|-|X_2||Y_2|+e(X_2,Y_2)+2\sqrt{tn}|B|\\[5pt]
&\leq \frac{(n-|B|)^2}{4}-\left(\frac{n}{2}-\sqrt{3tn}-2\right)^2+(k-3)n+2\sqrt{tn}|B|.
\end{align*}
Let $f(x):=\frac{(n-x)^2}{4}+2\sqrt{tn}x$. Recall that $f(x)$ is decreasing on $[1,n-4\sqrt{tn}]$.
 Since $n\geq 252t^2k>25t$ implies $ |B| \leq t+1\leq n-4\sqrt{tn}$, $f(|B|)\leq f(1) =\frac{(n-1)^2}{4}+2\sqrt{tn}$.
 Thus,
\begin{align*}
e(G)&\leq f(|B|)-\left(\frac{n}{2}-\sqrt{3tn}-2\right)^2+(k-3)n\\[5pt]
&\leq \frac{(n-1)^2}{4}+2\sqrt{tn} -\left(\frac{n}{2}-\sqrt{3tn}-2\right)^2+(k-3)n\\[5pt]
&\leq 2\sqrt{tn}+(\sqrt{3tn}+2)n-3tn-4\sqrt{3tn}+(k-3)n\\[5pt]
&\leq (k+\sqrt{3tn})n-\frac{t+1}{2}n.
\end{align*}
By $n\geq 252t^2k>20k+75t$ we have $(k+\sqrt{3tn})n\leq (\frac{n}{20}+\frac{n}{5})n= \frac{n^2}{4}$. It implies that $e(G)\leq \frac{n^2}{4}-\frac{t+1}{2}n <\frac{(n-t-1)^2}{4}$, a contradiction. Thus (iv) holds.

To show $(v)$, suppose for contradiction that there exist $uv\in E(B)$ and $x_1,x_2 \in X$ such that $ux_1,vx_2\in E(G)$. We distinguish two cases.

\vspace{5pt}
 {\bf \noindent Case 1. } $k=2$.

Then $|B|\leq t+1 \leq 2k-1=3$. By $u, v\in B$ we have $2\leq|B|\leq 3$.  Since $G$ is $C_5$-free, for every $ux',vx''\in E(G)$ with $x',x''\in X$ we have $N(x',Y)\cap N(x'',Y)=\emptyset$. It follows that
\begin{align}\label{ineq-3.4}
\deg(x',Y)+\deg(x'',Y)\leq |Y|.
\end{align}

If $|B|=3$, then $e(B)+e(B, X\cup Y)\leq2\sqrt{tn}|B|= 6\sqrt{tn}$. By \eqref{ineq-3.4}, $\deg(x_1,Y)+\deg(x_2,Y)\leq |Y|$. Since $|X|+|Y|-1=n-4$, it implies that
\[
e(X,Y) \leq|X\setminus \{x_1,x_2\}||Y|+\deg(x_1,Y)+\deg(x_2,Y) \leq(|X|-1)|Y|\leq\left\lfloor\frac{(n-4)^2}{4}\right\rfloor.
\]
Then
\begin{align*}
e(G)=e(B)+e(B, X\cup Y)+e(X,Y)&\leq6\sqrt{tn}+\left\lfloor\frac{(n-4)^2}{4}\right\rfloor\\[5pt]
&=\left\lfloor\frac{(n-3)^2}{4}\right\rfloor-\left\lfloor\frac{n-3}{2}\right\rfloor+6\sqrt{tn}.
\end{align*}
By $n\geq 252t^2k>144t$ we have $\left\lfloor\frac{n-3}{2}\right\rfloor+6>\frac{n}{2}\geq 6\sqrt{tn}$. It implies that $e(G)<\left\lfloor\frac{(n-3)^2}{4}\right\rfloor+6$,
a contradiction. Thus $|B|=2$, that is, $B=\{u,v\}$.

By (iv), $N(u, X\cup Y)\subset X$ and $N(v, X\cup Y)\subset X$. Applying \eqref{ineq-3.4} with $x'=x_1$ and $x''=x_2$,
\begin{align}\label{eq3.2}
\deg(x_1,Y)+\deg(x_2,Y) \leq |Y|.
\end{align}
If $\deg(u,X)+\deg(v,X)\leq 4$, then $e(B, X\cup Y)= \deg(u,X)+\deg(v,X)\leq 4$.
 As $|X|+|Y|-1=n-3$,
\[
\ e(X,Y)\leq |X\setminus \{x_1,x_2\}||Y|+\deg(x_1,Y)+\deg(x_2,Y)\leq (|X|-1)|Y|\leq \left\lfloor\frac{(n-3)^2}{4}\right\rfloor.
\]
Then
\begin{align*}
e(G)&=e(B)+e(B, X\cup Y)+e(X,Y)\leq 1+4+\left\lfloor\frac{(n-3)^2}{4}\right\rfloor< \left\lfloor\frac{(n-3)^2}{4}\right\rfloor+6,
\end{align*}
a contradiction. Thus we may assume that   $\deg(u,X)+\deg(v,X)\geq 5$.

Without loss of generality, we further assume that  $\deg(u,X)\geq \deg(v,X)$. Then $\deg(u,X)\geq 3$ and $e(B, X\cup Y)\leq 2\deg(u,X)$. Let $X_u=N(u,X)\setminus \{x_1, x_2\}$. Note that $x_2$ is a high-degree vertex. It follows that  $\deg(x_2,Y) \geq 2\sqrt{tn}-|\{u,v\}|=2\sqrt{tn}-2$. For every  $x\in X_u$, by applying \eqref{ineq-3.4} with $x'=x$ and $x''=x_2$ we get $\deg(x,Y)\leq |Y|-\deg(x_2,Y)\leq |Y|-2\sqrt{tn}+2$. By \eqref{eq3.2} we have
\begin{align*}
\ e(X,Y)&\leq |X\setminus (X_u\cup \{x_1,x_2\})||Y|+\deg(x_1,Y)+\deg(x_2,Y)+|X_u|(|Y|-2\sqrt{tn}+2)\\[5pt]
&\leq (|X|-1)|Y|-|X_u|(2\sqrt{tn}-2)\\[5pt]
&\leq \left\lfloor\frac{(n-3)^2}{4}\right\rfloor-(2\sqrt{tn}-2)(\deg(u,X)-2).
\end{align*}
Thus
\begin{align*}
e(G)&=e(B)+e(B, X\cup Y)+e(X,Y)\\[5pt]
&\leq 1+2\deg(u,X)+\left\lfloor\frac{(n-3)^2}{4}\right\rfloor-(2\sqrt{tn}-2)(\deg(u,X)-2)\\[5pt]
&\leq 1-(2\sqrt{tn}-4)\deg(u,X)+4\sqrt{tn}-4+\left\lfloor\frac{(n-3)^2}{4}\right\rfloor,
\end{align*}
Since $\deg(u,X)\geq 3$, we have
\[
e(G)\leq \left\lfloor\frac{(n-3)^2}{4}\right\rfloor-3(2\sqrt{tn}-4)+4\sqrt{tn}-3
=\left\lfloor\frac{(n-3)^2}{4}\right\rfloor-2\sqrt{tn}+9
< \left\lfloor\frac{(n-3)^2}{4}\right\rfloor+6,
\]
a contradiction.

\vspace{5pt}
 {\bf \noindent Case 2. } $k\geq3$.

Since $\deg(x_1), \deg(x_2)\geq 2\sqrt{tn}-t-1\geq \frac{3\sqrt{tn}}{2}+2$ and by  (iii), there exist distinct vertices $y_1\in N(x_1)\cap Y, y_2 \in N(x_2)\cap Y$ such that $\deg(y_1,X)\geq \frac{n}{2}-\sqrt{3tn}>\frac{n}{2}-2\sqrt{tn}$ and $\deg(y_2,X)\geq\frac{n}{2}-\sqrt{3tn}>\frac{n}{2}-2\sqrt{tn}$. Let $X_1=N(y_1,X)\cap N(y_2,X)\setminus \{x_1,x_2\}$, $Y_1=Y\setminus \{y_1,y_2\}$. Since $|X|\leq \frac{n}{2}+2\sqrt{tn}$ and $|Y|\geq \frac{n}{2}-2\sqrt{tn}$,
\begin{align*}
|X_1|&= |N(y_1,X)|+|N(y_2,X)|-|N(y_1,X)\cup N(y_2,X)|-2\\[5pt]
&\geq 2\left(\frac{n}{2}-2\sqrt{tn}\right)-|X|-2\\[5pt]
&\geq \frac{n}{2}-6\sqrt{tn}-2.
\end{align*}
Note also that $|Y_1|= |Y|-2\geq \frac{n}{2}-2\sqrt{tn}-2$. If there exists a path of length $2k-6$ with both ends in $X_1$, then together with $y_1x_1uvx_2y_2$ we obtain a cycle of length $2k+1$, a contradiction. Thus, there does not exist a path of length $2k-6$ with both ends in $X_1$. This implies that $G[X_1,Y_1]$ is $P_{2k-3}$-free. By Theorem \ref{erdos-gallai-1},
\[
e(X_1,Y_1)\leq \frac{(2k-3)-1}{2}(|X_1|+|Y_1|)\leq(k-2)n.
\]
Then
\begin{align*}
e(G)&\leq |X||Y|-|X_1||Y_1|+e(X_1,Y_1)+(t+1)\cdot 2\sqrt{tn}\\[5pt]
&\leq \frac{(n-2)^2}{4}-\left(\frac{n}{2}-6\sqrt{tn}-2\right)
\left(\frac{n}{2}-2\sqrt{tn}-2\right)+(k-2)n+2(t+1)\sqrt{tn}\\[5pt]
&\leq \frac{(n-2)^2}{4}-\frac{n^2}{4}+4n\sqrt{tn}-12tn+2n+(k-2)n+2(t+1)\sqrt{tn}\\[5pt]
&\leq (4\sqrt{tn}+k)n-12tn+2(t+1)\sqrt{tn}\\[5pt]
&\leq (4\sqrt{tn}+k)n-\frac{t+1}{2}n-10tn+4t\sqrt{tn}.
\end{align*}
It is easy to verify that $10tn\geq 4t\sqrt{tn}$. By $n\geq 252t^2k>20k+400t$ we have $(4\sqrt{tn}+k)n\leq(\frac{n}{5}+\frac{n}{20})n=\frac{n^2}{4}$. It implies that $e(G)\leq \frac{n^2}{4}-\frac{t+1}{2}n < \frac{(n-t-1)^2}{4}$, a contradiction. Thus (v) holds.

To show $(vi)$, let $uwv$ be a path in $B$ and $\deg(w, X\cup Y)\leq 1$, assume that there exist $x\in X, y\in Y$ such that $uy,vx\in E(G)$. We distinguish two cases.

\vspace{5pt}
 {\bf \noindent Case 1. } $k=2$.

 Since $|B|\leq t+1\leq 2k-1=3$, $B=\{u,v,w\}$. By (iv), we have $N(u, X\cup Y)\subset Y$ and $N(v, X\cup Y)\subset X$. It follows that
\[
e(B, X\cup Y)\leq \deg(u,Y)+\deg(v,{X})+\deg(w, X\cup Y)\leq \deg(u,{Y})+\deg(v,{X})+1.
\]
Since $G$ is $C_5$-free, there is no edge between $N(u,Y)$ and $N(v,X)$. That is,
 \begin{align*}
e(X,Y)\leq |X||Y|&-\deg(u,{Y})\deg(v,{X}).
\end{align*}
It follows that
\begin{align*}
e(G)&=e(B)+e(B, X\cup Y)+e(X,Y)\\[5pt]
&\leq 3+\deg(u,{Y})+\deg(v,{X})+1+|X||Y|-\deg(u,{Y})\deg(v,{X})\\[5pt]
&= 5+|X||Y|-(\deg(u,{Y})-1)(\deg(v,{X})-1).
\end{align*}
Since $|X|+|Y|= n-3$, $|X||Y|\leq \left\lfloor\frac{(n-3)^2}{4}\right\rfloor$. Moreover, $\deg(u,{Y})\geq 1$ and $\deg(v,{X})\geq 1$, we have $e(G)\leq 5+\left\lfloor\frac{(n-3)^2}{4}\right\rfloor$. It follows that
$e(G)<\left\lfloor\frac{(n-2)^2}{4}\right\rfloor+3$ for $t=1$ and $e(G)< \left\lfloor\frac{(n-3)^2}{4}\right\rfloor+6$ for $t=2$, a contradiction.

\vspace{5pt}
 {\bf \noindent Case 2. } $k\geq3$.

 Let $X_1=N(y, X)\setminus \{x\}$ and $Y_1=N(x,Y)\setminus \{y\}$. Recall that $x,y$ are both high-degree vertices. Since $n\geq 252t^2k>  64t$ implies $|B|\leq t+1 \leq 2t \leq \frac{\sqrt{tn}}{2}-2$, we infer that
 \[
   |X_1|, |Y_1|\geq 2\sqrt{tn}-|B|-1\geq \frac{3\sqrt{tn}}{2}+1>\sqrt{tn}.
 \]
It follows that  $|X_1| |Y_1|> tn$.

If $k=3$ then  by the $C_7$-free property there is no edge between $X_1$ and $Y_1$.
Thus,
\begin{align*}
e(G)&=e(B)+e(B, X\cup Y)+e(X,Y)\\[5pt]
&\leq(t+1)\cdot 2\sqrt{tn}+ |X||Y|-|X_1||Y_1|\\[5pt]
&\leq \frac{(n-3)^2}{4}-tn+2(t+1)\sqrt{tn}\\[5pt]
&\leq \frac{(n-t-1)^2}{4}-\frac{t+1}{2}n+2(t+1)\sqrt{tn}.
\end{align*}
By $n\geq 252t^2k\geq 16t$ we have $\frac{t+1}{2}n\geq 2(t+1)\sqrt{tn}$. Thus $e(G)\leq \frac{(n-t-1)^2}{4}$, a contradiction.

If $k\geq 4$, by (iii) there exist  $x_1 \in X_1$,  $y_1\in Y_1$ such that $\deg(x_1,{Y}) >\frac{n}{2}-2\sqrt{tn}$ and $\deg(y_1,{X}) >\frac{n}{2}-2\sqrt{tn}$. Let $X_2=N(y_1, X)\setminus \{x,x_1\}$ and $Y_2=N(x_1, Y)\setminus \{y,y_1\}$. By the $C_{2k+1}$-free property there is no a path of length $2k-7$ with one end in $X_2$ and the other one in $Y_2$. By Theorem \ref{erdos-gallai-1},
\[
e(X_2,Y_2)\leq \frac{(2k-7)-1}{2}(|X_1|+|Y_1|)\leq(k-4)n.
\]
Using $|X_2|, |Y_2|\geq \frac{n}{2}-2\sqrt{tn}-2$, we have
\begin{align*}
e(G)&=e(B)+e(B, X\cup Y)+e(X,Y)\\[5pt]
&\leq(t+1)\cdot 2\sqrt{tn}+ |X||Y|-|X_2||Y_2|+e(X_2,Y_2)\\[5pt]
&\leq \left\lfloor\frac{(n-3)^2}{4}\right\rfloor-\left(\frac{n}{2}-2\sqrt{tn}-2\right)^2+(k-4)n+2(t+1)\sqrt{tn}\\[5pt]
&\leq (2\sqrt{tn}+2)n-4tn+(k-4)n+4t\sqrt{tn}\\[5pt]
&\leq (k+2\sqrt{tn})n-\frac{t+1}{2}n-2tn+4t\sqrt{tn}.
\end{align*}
It is easy to verify that $2tn\geq 4t\sqrt{tn}$. By $n\geq 252t^2k>10k+64t$ we have $(k+2\sqrt{tn})n\leq (\frac{n}{10}+\frac{n}{8})n< \frac{n^2}{4}$. It implies that $e(G)\leq \frac{n^2}{4}-\frac{t+1}{2}n <\frac{(n-t-1)^2}{4}$, a contradiction. Thus (vi) holds and the lemma is proven.
\end{proof}

\section{The Proof of Theorem \ref{thm-main1}}

\begin{proof}[Proof of Theorem \ref{thm-main1}]
Let $G$ be a $C_{2k+1}$-free graph on $n$ vertices with $d_2(G)\geq t$. If $e(G)< \lfloor\frac{(n-t-1)^2}{4}\rfloor+\binom{t+2}{2}$ then we have nothing to do. Thus we assume that
\begin{align}\label{ineq-asmption}
e(G)\geq \frac{(n-t-1)^2}{4}+\binom{t+2}{2}
\end{align}
and we are left to show that $G=H(n,t)$ up to isomorphism.
By Lemma \ref{lem-2.2}, there exists $B\subseteq V(G)$ with $t\leq|B| \leq t+1$ such that $G-B$ is a bipartite graph on partite sets $X,Y$ with (i)-(vi) of Lemma \ref{lem-2.2}  hold.

\begin{claim}\label{claim_1}
Let $B'\subseteq  B$. If  $G-B'$ is non-bipartite, then each shortest odd cycle  in $G-B'$ has length at most $2k-1$.
\end{claim}
\begin{proof}
Since $G-B$ is bipartite and  $G-B'$ is non-bipartite, $|B'|\leq |B|-1\leq t$. Let $C$ be a shortest cycle in $G-B'$, let $S=V(C)$. Suppose for contradiction  that  $|S|=2\ell+1\geq 2k+3$. By Fact \ref{fact-3} each $v\in B'$ has at most $\ell$ neighbors on $C$.  By Fact \ref{fact-1} each $v\in V(G)\setminus (B'\cup S)$ has at most two neighbors on $C$. It follows that
  \begin{align*}
  e(S)+e(G-S,S)&\leq 2\ell+1+ {2\left(n-|B'|-|S|\right)+|B'| \ell}\\[5pt]
  &\leq 2\ell+1+ {2\left(n-t-2\ell-1\right)+t \ell}\\[5pt]
  &=2n +(t-2)\ell-2t-1.
  \end{align*}
By Fact \ref{fact-2} we have $2\ell+1\leq \frac{2n}{3}$. It follows that $n-2\ell-1\geq \frac{n}{3} \geq 4k$. Since $G-S$ is $C_{2k+1}$-free, by Theorem \ref{thm-1.1} we have  $e(G-S)\leq{\frac{(n-2\ell-1)^2}{4}}$.
Thus,
\begin{align*}
e(G)&=e(S)+e(G-S,S)+e(G-S)\leq 2n +(t-2)\ell-2t-1+\frac{(n-2\ell-1)^2}{4}=:f(\ell).
\end{align*}
Since $f(\ell)$ is a convex function of $\ell$ with  $k+1\leq \ell\leq \frac{n-1}{2}$, $e(G)\leq \max\{f(\frac{n-1}{2}),f(k+1)\}$.
By $n\geq 318t^2k>14t+2$ we have $\frac{n-2(t+1)}{4}\geq 3t$. It follows that
\[
f\left(\frac{n-1}{2}\right) =2n +\frac{(t-2)(n-1)}{2}-2t-1\leq 3tn\leq \frac{(n-2(t+1))n}{4}< \frac{(n-t-1)^2}{4}.
\]
For $t\leq 2k-3$,
\begin{align*}
f(k+1) &=2n +(t-2)(k+1)-2t-1+\frac{(n-2k-3)^2}{4}\\[5pt]
&=2n+(t-2)(k+1)-2t-1+\frac{(n-2k+2)^2}{4}-\frac{5}{2}n+\frac{(2k+3)^2-(2k-2)^2}{4}\\[5pt]
&< \frac{(n-2k+2)^2}{4}-\frac{n}{2}+(t-2)k-2+\frac{5(4k+1)}{4} \\[5pt]
&< \frac{(n-2k+2)^2}{4}-\frac{n}{2}+(t+3)k.
\end{align*}
Since $n\geq 318t^2k$ implies $\frac{n}{2}\geq (t+3)k$, we infer that  $f(k+1)< \frac{(n-t-1)^2}{4}$. For $t=2k-2$,
\begin{align*}
f(k+1) &= 2n+2(k-2)(k+1)-4k+3+\frac{(n-2k-3)^2}{4}\\[5pt]
&=\frac{(n-2k+1)^2}{4}+2k^2-2k+1\\[5pt]
&\leq\frac{(n-2k+1)^2}{4} +\binom{2k}{2}-1\\[5pt]
&= \frac{(n-t-1)^2}{4}+\binom{t+2}{2}-1.
\end{align*}
Thus $e(G)\leq \max\{f(\frac{n-1}{2}),f(k+1)\}\leq\frac{(n-t-1)^2}{4}+\binom{t+2}{2}-1$, contradicting \eqref{ineq-asmption}.
\end{proof}
\begin{claim}\label{claim_2}
Let $B'\subseteq  B$ such that $G-B'$ is non-bipartite and let $C$ be a shortest  odd cycle in $G-B'$. Then
there is at most one high-degree vertex on $C$.
\end{claim}
\begin{proof}
Suppose that $C=v_0v_1\ldots v_{2\ell}v_0$ and there are two high-degree vertices  on $C$, say $v_i, v_j$ with $0\leq i<j\leq 2\ell$. Set $Q_1=v_iv_{i+1}\ldots v_j$  and $Q_2=v_jv_{j+1}\ldots v_{2\ell} v_0\ldots v_i$. Let $|V(Q_1)|=\ell_1+1$ and  $|V(Q_2)|=\ell_2+1$. By Claim \ref{claim_1} we know that $\ell\leq k-1$. It follows that  $\ell_1,\ell_2\in \{1,2,\ldots,2k-2\}$. Without loss of generality, assume that $\ell_1$ is odd and $\ell_2$ is even.

Note that $C$ is an induced cycle and each $v_i$ has at most $2$ neighbors on $C$.  Let $X_1=X\setminus S$ and $Y_1=Y\setminus S$.
Since $n\geq 318t^2k>144t$ implies  $\deg(v_i)\geq 2\sqrt{tn} \geq \frac{3\sqrt{tn}}{2}+ 6t$,  we have
\[
\deg(v_i, X_1\cup Y_1)\geq \deg(v_i) -|B| -\deg(v_i,S)  \geq \frac{3\sqrt{tn}}{2}+ 6t -(t+1)-2\geq \frac{3\sqrt{tn}}{2}+2.
\]
Similarly, $\deg(v_j, X_1\cup Y_1)  \geq \frac{3\sqrt{tn}}{2}+2$. By  Lemma \ref{lem-2.2} (iii), there exist distinct vertices $x_1\in N(v_i,X_1\cup Y_1), x_2 \in N(v_j,X_1\cup Y_1)$ such that $\deg(x_1,X\cup Y)> \frac{n}{2}-2\sqrt{tn}$ and $\deg(x_2,X\cup Y)> \frac{n}{2}-2\sqrt{tn}$. Then by Claim \ref{claim_1}
 \[
 \deg(x_1,X_1\cup Y_1)\geq \deg(x_1,X\cup Y)-|S|\geq \frac{n}{2}-2\sqrt{tn}-2k+1.
 \]
 Similarly, $\deg(x_2,X_1\cup Y_1)\geq \frac{n}{2}-2\sqrt{tn}-2k+1$. Now we distinguish two cases.

\vspace{5pt}
 {\bf \noindent Case 1. } $x_1$ and $x_2$ belong to the same one of $X_1,Y_1$.

Without loss of generality, assume $x_1, x_2\in X_1$. Let $Y_2=N(x_1,Y_1)\cap N(x_2,Y_1)$ and $X_2=X_1\setminus\{x_1,x_2\}$. By Lemma \ref{lem-2.2} (ii)  $|X_2|\geq |X|-|S|-2\geq \frac{n}{2}-2\sqrt{tn}-2k-1$. By Lemma \ref{lem-2.2} (ii) $|Y_1|\leq |Y|\leq \frac{n}{2}+2\sqrt{tn}$. Then
\begin{align*}
|Y_2|=|N(x_1,Y_1)\cap N(x_2,Y_1)|&\geq \deg(x_1,Y_1)+\deg(x_2,Y_1)-|Y_1|\\[5pt]
&\geq 2\left(\frac{n}{2}-2\sqrt{tn}-2k\right)-|Y_1|\\[5pt]
&\geq \frac{n}{2}-6\sqrt{tn}-4k.
\end{align*}

Let
\[
L=\left\{v\in V(G)\setminus B\colon \deg(v,X\cup Y)\leq \frac{n}{2}-\sqrt{3tn}\right\}.
\]
By Lemma \ref{lem-2.2} (iii) we have $|L| \leq \frac{3\sqrt{tn}}{2}$. It follows that   $G[X\setminus L,Y\setminus L]$ is a subgraph with minimum degree at least $\frac{n}{2}-\sqrt{3tn}-|L|>\frac{n}{2}-4\sqrt{tn}$. Note that
\[
|Y\setminus Y_2|\leq \left(\frac{n}{2}+2\sqrt{tn}\right) -\left(\frac{n}{2}-6\sqrt{tn}-4k\right) =8\sqrt{tn}+4k.
\]
Then $G':=G[X_2\setminus L,Y_2\setminus L]$ is a graph on $|X_2|+|Y_2|-|L|$ vertices with
\[
\delta(G')\geq \frac{n}{2}-4\sqrt{tn}-|Y\setminus Y_2|\geq \frac{n}{2}-12\sqrt{tn}-4k.
 \]
Since $n\geq 318t^2k\geq \frac{318^2}{153^2}\times 12^2t$ implies $\sqrt{tn} \leq \frac{1}{12} \times \frac{153}{318}n$, we have $\frac{n}{2}-12\sqrt{tn}\geq \frac{6}{318}n\geq 6k$. It follows that $\delta(G')\geq 2k$. Thus there exists a path of length at least $\delta(G')+1\geq 2k+1$. For $\ell_1\leq 2k-4$,  we  choose a sub-path $Q$ of length $2k-3-\ell_1$ with both ends in $Y_2$. For $\ell_1=2k-3$, we choose $Q$ as an arbitrary vertex in $Y_1$. Then $v_jx_2Qx_1Q_1$ forms a cycle of length $2k+1$, a contradiction.

\vspace{5pt}
 {\bf \noindent Case 2. } $x_1$ and $x_2$ belong to different ones of $X_1$, $Y_1$.

 Without loss of generality, assume that $x_1\in X_1$ and $x_2\in Y_1$. Let $X_2=N(x_2,X_1)\setminus \{x_1\}$ and $Y_2=N(x_1,Y_1)\setminus \{x_2\}$. Then $|X_2|\geq \deg(x_2,X_1)-1\geq \frac{n}{2}-2\sqrt{tn}-2k$. Similarly, $|Y_2|\geq \frac{n}{2}-2\sqrt{tn}-2k$. Let
\[
L=\left\{v\in V(G)\setminus B\colon \deg(v,X\cup Y)\leq \frac{n}{2}-\sqrt{3tn}\right\}.
\]By Lemma \ref{lem-2.2} (iii) we have $|L| \leq \frac{3\sqrt{tn}}{2}$. It follows that   $G[X\setminus L,Y\setminus L]$ is a subgraph with minimum degree at least $\frac{n}{2}-\sqrt{3tn}-|L|>\frac{n}{2}-4\sqrt{tn}$. Note that
\[
|Y\setminus Y_2|\leq \left(\frac{n}{2}+2\sqrt{tn}\right) -\left(\frac{n}{2}-2\sqrt{tn}-2k\right) =4\sqrt{tn}+2k.
\]
Then $G':=G[X_2\setminus L,Y_2\setminus L]$ is a graph on $|X_2|+|Y_2|-|L|$ vertices with \[
\delta(G')\geq \frac{n}{2}-4\sqrt{tn}-|Y\setminus Y_2|\geq \frac{n}{2}-8\sqrt{tn}-2k.
 \]
Since $n\geq 318t^2k>\frac{318^2}{155^2}\times 64t$ implies $\sqrt{tn} \leq \frac{1}{8} \times \frac{155}{318}n$, we have $\frac{n}{2}-8\sqrt{tn}\geq \frac{4}{318}n\geq 4k$. It implies that $\delta(G')\geq 2k$. Then there exists a path of length at least $\delta(G')+1\geq 2k+1$. For $\ell_2\leq 2k-4$, we choose a sub-path $Q$ of length $2k-3-\ell_2$ with one end in $X_1$ and the other one in $Y_1$.  Then $v_ix_1Qx_2Q_2$ forms a cycle of length $2k+1$, a contradiction.

We are left with the case $\ell_2=2k-2$. The $C_{2k+1}$-free property implies that there is no edge between $N(v_i, X_1)$ and $N(v_j,  Y_1)$. Since $n\geq 318t^2k\geq \frac{16}{(2-\frac{23}{25})^2}t$ implies $4t\leq \left(2-\frac{23}{25}\right)\sqrt{tn}$,
\[
\deg(v_i, X_1)\geq \deg(v_i) -|B| -\deg(v_i,S)\geq 2\sqrt{tn} -t-3\geq2\sqrt{tn} -4t\geq \frac{23}{25}\sqrt{tn}.
\]
Similarly, $\deg(v_j, Y_1)\geq \frac{23}{25}\sqrt{tn}$.
Then
\[
e(G-B)\leq |X||Y| -|N(v_i, X_1)||N(v_j, Y_1)|\leq \frac{(n-|B|)^2}{4}-\frac{5tn}{6}.
\]

By Lemma \ref{lem-2.2} we have
\begin{align*}
e(G)=e(G-B)+ 2\sqrt{tn}|B|\leq \frac{(n-|B|)^2}{4} -\frac{5tn}{6}+2\sqrt{tn}|B|.
\end{align*}
Let $f(x):=\frac{(n-x)^2}{4}+2\sqrt{tn}x$. It is easy to verify that $f(x)$ is decreasing on $[1,n-4\sqrt{tn}]$.  Note that $n\geq 318t^2k>25t$ implies $ |B| \leq t+1 \leq n-4\sqrt{tn}$. By $|B|\geq t$ we infer that
\[
e(G) \leq f(t)-\frac{5tn}{6}=\frac{(n-t)^2}{4}+2t\sqrt{tn}-\frac{5tn}{6}.
\]
Since $n\geq 36t$ we have $2t\sqrt{tn}\leq \frac{tn}{3}$. Then $e(G)\leq\frac{(n-t)^2}{4}-\frac{tn}{2}
\leq\frac{(n-t-1)^2}{4}+\frac{n-t}{2}-\frac{tn}{2}\leq\frac{(n-t-1)^2}{4}$, contradicting \eqref{ineq-asmption}. Thus there is an edge $xy$ between $N(v_i, X_1\cup Y_1)$ and $N(v_j, X_1\cup Y_1)$. It follows that
 $xyQ_2x$ is a cycle of length $2k+1$, a contradiction. Thus the claim holds.
\end{proof}

Now we choose $B_0\subset B$ such that $G-B_0$ is bipartite and $|B_0|$ is minimal.  By $d_2(G)\geq t$ we have $t \leq |B_0|\leq t+1$. By the minimality of $|B_0|$, we infer that for each $u\in B_0$, $G-(B_0\setminus \{u\})$ is not bipartite.  Let $C_{u}$ be a shortest odd cycle in $G-(B_0\setminus \{u\})$. By Claim \ref{claim_2}, there is at most one high-degree vertex on $C_{u}$. Since $|V(C_u)|\geq 3$, $C_{u}$ contains at least two low-degree vertices. As $V(C_u)\cap B_0=\{u\}$, we infer that $C_{u}$ contains a low-degree vertex that is not in $B_0$, implying that $B\setminus B_0\neq \emptyset$. By $t\leq |B_0|\leq |B|\leq t+1$ we have $|B|= t+1$ and  $|B_0|=t$. Then $C_{u}$ contains exactly two low-degree vertices and one high-degree vertex, that is, $|V(C_{u})|=3$.  Let $B\setminus B_0=\{u_0\}$ and let $x_u$ be a high-degree vertex on $C_u$.  Then $V(C_u)=\{u,u_0,x_u\}$. Since $u$ is chosen from $B_0$ arbitrarily,  we conclude that  $u_0u\in E(G)$ for all $u\in B_0$.

Assume that $B=\{u_0,u_1,u_2,\ldots,u_t\}$ and let $B_i=B\setminus \{u_i\}$, $i=1,2,\ldots,t$.  We claim that $G-B_i$ is bipartite. Indeed, otherwise let $C^{i}$ be a shortest odd cycle in $G-B_i$. Note that $|B_i|= t$. By applying Claim \ref{claim_2}  to $C^i$, $C^i$ contains at most one high-degree vertex. It follows that $C^i$ contains at least two low-degree vertices. But there is only one low-degree vertex $u_i$ in  $G-B_i$, a contradiction. Thus $G-B_i$ is bipartite. Let $B_{ij} =B\setminus \{u_i,u_j\}$ for $0\leq i<j\leq t$.  By $d_2(G)\geq t$ and $|B_{ij}|=t-1$, we see that $G-B_{ij}$ is non-bipartite.

%Note that $V(C^i)\cap B=\{u_i\}$. It implies that $C^i$ contains at lease one low-degree vertex that are not in $B$, contradicting the definition of $B$. Thus $G-B_i$ is bipartite. By $|B_i|=t$ and $d_2(G)\geq t$, $|B_i|$ is minimal such that $G-B_i$ is bipartite for  $i=1,2,\ldots,t$. Then $G-(B_i\setminus u_j)$ is not bipartite for each $u_j\in B_i$.  Let $C_{u_j}$ be a shortest odd cycle in $G-(B_i\setminus u_j)$.
% Now apply  Claim \ref{claim_2}  to $C_{u_j}$, we obtain that $C_{u_j}$ contains at least two low-degree vertices. By the same argument, we have  $u_iu_j\in E(G)$ for all $j\in \{0,1,\ldots,t\}\setminus \{i\}$.

\begin{claim}\label{claim_3}
For $0\leq i<j\leq t$,
there exists a high-degree vertex $x_{ij}$ such that $u_iu_jx_{ij}$ forms a triangle in $G$.
\end{claim}
\begin{proof}
Let $C$ be a shortest odd cycle in $G-B_{ij}$.  By Claim \ref{claim_2}, $C$ contains  at most one high-degree vertex. Since $G-B_{ij}$ contains exactly two low-degree vertices $u_i,u_j$, $C$ has to be a triangle on $\{u_i,u_j,x_{ij}\}$ for some high-degree vertex $x_{ij}$.
\end{proof}

By Claim \ref{claim_3}, $B$ spans a clique of size $t+1$ in $G$. Moreover, for each $u_iu_j\in E(B)$ there exists a high-degree vertex $x_{ij}$ such that  $C_{ij}=u_iu_jx_{ij}u_i$ is a triangle in $G$.

\begin{claim}\label{claim_4}
For each $u_iu_j\in E(B)$,  $N(u_i,X\cup Y)=N(u_j,X\cup Y)=\{x_{ij}\}$.
\end{claim}
\begin{proof}
Otherwise  there exists a high-degree vertex $y\neq x_{ij}$ such that $u_ix_{ij},u_jy\in E(G)$. By  Lemma \ref{lem-2.2} (iv), it implies that $x_{ij},y$  belong to same side of $X$ or $Y$. By symmetry assume that $y, x_{ij}\in X$, contradicting (v) of Lemma \ref{lem-2.2}.
\end{proof}

Recall that $G[B]$ is a clique of size $t+1$. By Claim \ref{claim_4}, $x_{ij}=x_{i'j'}=:x^*$ for all $u_iu_j,u_{i'}u_{j'}\in E(B)$. That is, each vertex in $B$ is adjacent to the unique high-degree vertex $x^*$.
Hence $G$ is isomorphic to a subgraph of $H(n,t)$. By \eqref{ineq-asmption} we conclude that $G=H(n,t)$ up to isomorphism
 and the theorem is proven.
\end{proof}

\section{The Proof of Theorem  \ref{thm-main2}}

\begin{proof}[Proof of Theorem  \ref{thm-main2}]
Let $G$ be a $C_{2k+1}$-free on $n$ vertices with
\begin{align}\label{ineq-asmption2}
\gamma_2(G)\geq \binom{\lfloor\frac{t+2}{2}\rfloor}{2}+\binom{\lceil\frac{t+2}{2}\rceil}{2}
\overset{\eqref{equalities2.6}}{=}f(t+1).
\end{align}
If $e(G)<  \left\lfloor\frac{(n-t-1)^2}{4}\right\rfloor+\binom{t+2}{2}$ then we are done. Thus we may assume that
\begin{align}\label{ineq-asmption1}
e(G)\geq  \left\lfloor\frac{(n-t-1)^2}{4}\right\rfloor+\binom{t+2}{2}
\end{align}
and we are left to show that $G=H(n,t)$ up to isomorphism.
By Lemma \ref{lem-2.2}, there exists $B\subseteq V(G)$ with $|B| \leq t+1$ such that $G-B$ is a bipartite graph on partite sets $X, Y$ with (i)-(vi) of Lemma \ref{lem-2.2}  hold. Let
\begin{align*}
&B_{11}=\{u\in B\colon \deg(u,Y)= 1\}, \quad B_{12}=\{u\in B\colon \deg(u,X)= 1\}, \\[5pt]
&X_{2}=\{u\in B\colon \deg(u,Y)\geq 2\}, \quad Y_{2}=\{u\in B\colon \deg(u,X)\geq 2\}, \\[5pt]
&B_0=\{u\in B\colon \deg(u, X\cup Y)=0\}.
\end{align*}
By Lemma \ref{lem-2.2} (iv), $B_{11}\cap B_{12}=\emptyset$, $B_{11}\cap Y_2=\emptyset$, $X_2\cap B_{12}=\emptyset$ and  $X_2\cap Y_2=\emptyset$. It follows that $(B_0, B_{11}, B_{12}, X_{2}, Y_{2})$ is partition of $B$. Furthermore, we  have
\begin{align}\label{eq_5.2}
e(X_2,X)=0, \ e(Y_2, Y)=0,\ e(B_{11},X)=0 ,\ e(B_{12},Y)=0.
\end{align}

\begin{claim}\label{claim-5}
$X_2$ and $Y_2$ are both independent sets of $G$.
\end{claim}
\begin{proof}
Suppose not, by symmetry assume that $uv\in E(X_2)$. Since $\deg(u,Y)\geq 2$ and $\deg(v,Y)\geq 2$, there are different vertices $y_1,y_2 \in Y$ such that $uy_1,vy_2\in E(G)$, contradicting (v) of Lemma \ref{lem-2.2}.  Thus $X_2$ is an independent set. Similarly, $Y_2$ is an independent set.
\end{proof}
\begin{claim}\label{claim-6}
$e(B_{11},X_2)=0$ and $e(B_{12},Y_2)=0$.
\end{claim}
\begin{proof}
Suppose that there exist $u\in B_{11}$, $v\in X_2$ such that $uv\in E(G)$. Since $\deg(u,Y)= 1$ and $\deg(v,Y)\geq 2$, one can find different vertices $y_1, y_2\in Y$ such that $uy_1, vy_2 \in E(G)$, contradicting Lemma \ref{lem-2.2} (v). Thus $e(B_{11},X_2)=0$. Similarly, $e(B_{12},Y_2)=0$.
\end{proof}

Set $|B_0|=z$, $|B_{11}|=p$, $|B_{12}|=q$, $|X_2|=x$ and $|Y_2|=y$. Clearly $x+y+z+p+q=|B|\leq t+1$.  Now we further partition  $B_0$ to $X_0$ and $Y_0$ such that $|X_0|=\left\lfloor\frac{z}{2}\right\rfloor$ and $|Y_0|=\left\lceil\frac{z}{2}\right\rceil$. Partition $B_{11}$ to $X_{11}$, $Y_{11}$ such that $|X_{11}|=\left\lceil\frac{p}{2}\right\rceil$, $|Y_{11}|=\left\lfloor\frac{p}{2}\right\rfloor$ and partition $B_{12}$ to $X_{12}$, $Y_{12}$ such that $|X_{12}|=\left\lfloor\frac{q}{2}\right\rfloor$, $|Y_{12}|=\left\lceil\frac{q}{2}\right\rceil$.
Let
\[
\widetilde{X}=X_0\cup X_{11}\cup X_{12}\cup X_2\cup X,\  \widetilde{Y}=Y_0\cup Y_{11}\cup Y_{12}\cup Y_2\cup Y.
\]
Clearly, $(\widetilde{X},\widetilde{Y})$ is a bi-partition of $V(G)$ (as shown in Figure \ref{fig-1}) and $\gamma_2(G)\leq e(\widetilde{X})+e(\widetilde{Y})$.

\begin{figure}[H]
\centering
\ifpdf
  \setlength{\unitlength}{0.05 mm}%
  \begin{picture}(1600.4, 851.8)(0,0)
  \put(0,0){\includegraphics{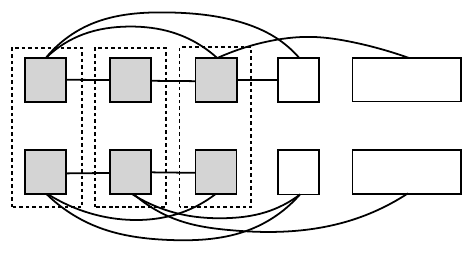}}
  \put(114.37,703.39){\fontsize{7.11}{8.54}\selectfont \makebox(75.0, 50.0)[l]{$B_0$\strut}}
  \put(381.42,705.08){\fontsize{7.11}{8.54}\selectfont \makebox(100.0, 50.0)[l]{$B_{11}$\strut}}
  \put(678.89,705.08){\fontsize{7.11}{8.54}\selectfont \makebox(100.0, 50.0)[l]{$B_{12}$\strut}}
  \put(970.17,239.44){\fontsize{7.11}{8.54}\selectfont \makebox(75.0, 50.0)[l]{$Y_2$\strut}}
  \put(93.98,443.71){\fontsize{7.11}{8.54}\selectfont \makebox(75.0, 50.0)[l]{ $\left\lfloor\frac{z}{2}\right\rfloor$\strut}}
  \put(115.84,367.21){\fontsize{7.11}{8.54}\selectfont \makebox(75.0, 50.0)[l]{$\left\lceil\frac{z}{2}\right \rceil $\strut}}
  \put(395.13,443.71){\fontsize{7.11}{8.54}\selectfont \makebox(75.0, 50.0)[l]{$\left\lceil\frac{p}{2}\right \rceil $\strut}}
  \put(397.56,367.21){\fontsize{7.11}{8.54}\selectfont \makebox(75.0, 50.0)[l]{$\left\lfloor\frac{p}{2}\right\rfloor$\strut}}
  \put(685.35,441.29){\fontsize{7.11}{8.54}\selectfont \makebox(75.0, 50.0)[l]{$\left\lfloor\frac{q}{2}\right\rfloor$\strut}}
  \put(686.57,367.21){\fontsize{7.11}{8.54}\selectfont \makebox(75.0, 50.0)[l]{$\left\lceil\frac{q}{2}\right \rceil$\strut}}
  \put(987.72,438.86){\fontsize{7.11}{8.54}\selectfont \makebox(25.0, 50.0)[l]{$x$\strut}}
  \put(992.57,364.78){\fontsize{7.11}{8.54}\selectfont \makebox(25.0, 50.0)[l]{$y$\strut}}
  \put(974.52,557.09){\fontsize{7.11}{8.54}\selectfont \makebox(75.0, 50.0)[l]{$X_2$\strut}}
  \put(1371.33,549.41){\fontsize{7.11}{8.54}\selectfont \makebox(25.0, 50.0)[l]{$X$\strut}}
  \put(1368.75,239.44){\fontsize{7.11}{8.54}\selectfont \makebox(25.0, 50.0)[l]{$Y$\strut}}
  \put(113.57,555.32){\fontsize{7.11}{8.54}\selectfont \makebox(75.0, 50.0)[l]{$X_0$\strut}}
  \put(392.28,562.47){\fontsize{7.11}{8.54}\selectfont \makebox(100.0, 50.0)[l]{$X_{11}$\strut}}
  \put(678.56,562.47){\fontsize{7.11}{8.54}\selectfont \makebox(100.0, 50.0)[l]{$X_{12}$\strut}}
  \put(395.46,243.06){\fontsize{7.11}{8.54}\selectfont \makebox(100.0, 50.0)[l]{$Y_{11}$\strut}}
  \put(686.81,239.44){\fontsize{7.11}{8.54}\selectfont \makebox(100.0, 50.0)[l]{$Y_{12}$\strut}}
  \put(118.75,250.60){\fontsize{7.11}{8.54}\selectfont \makebox(75.0, 50.0)[l]{$Y_0$\strut}}
  \end{picture}%
\else
  \setlength{\unitlength}{0.05 mm}%
  \begin{picture}(1600.4, 851.8)(0,0)
  \put(0,0){\includegraphics{6}}
  \put(114.37,703.39){\fontsize{7.11}{8.54}\selectfont \makebox(75.0, 50.0)[l]{$B_0$\strut}}
  \put(381.42,705.08){\fontsize{7.11}{8.54}\selectfont \makebox(100.0, 50.0)[l]{$B_{11}$\strut}}
  \put(678.89,705.08){\fontsize{7.11}{8.54}\selectfont \makebox(100.0, 50.0)[l]{$B_{12}$\strut}}
  \put(970.17,239.44){\fontsize{7.11}{8.54}\selectfont \makebox(75.0, 50.0)[l]{$Y_2$\strut}}
  \put(93.98,443.71){\fontsize{7.11}{8.54}\selectfont \makebox(75.0, 50.0)[l]{ $\left\lfloor\frac{z}{2}\right\rfloor$\strut}}
  \put(115.84,367.21){\fontsize{7.11}{8.54}\selectfont \makebox(75.0, 50.0)[l]{$\left\lceil\frac{z}{2}\right \rceil $\strut}}
  \put(395.13,443.71){\fontsize{7.11}{8.54}\selectfont \makebox(75.0, 50.0)[l]{$\left\lceil\frac{p}{2}\right \rceil $\strut}}
  \put(397.56,367.21){\fontsize{7.11}{8.54}\selectfont \makebox(75.0, 50.0)[l]{$\left\lfloor\frac{p}{2}\right\rfloor$\strut}}
  \put(685.35,441.29){\fontsize{7.11}{8.54}\selectfont \makebox(75.0, 50.0)[l]{$\left\lfloor\frac{q}{2}\right\rfloor$\strut}}
  \put(686.57,367.21){\fontsize{7.11}{8.54}\selectfont \makebox(75.0, 50.0)[l]{$\left\lceil\frac{q}{2}\right \rceil$\strut}}
  \put(987.72,438.86){\fontsize{7.11}{8.54}\selectfont \makebox(25.0, 50.0)[l]{$x$\strut}}
  \put(992.57,364.78){\fontsize{7.11}{8.54}\selectfont \makebox(25.0, 50.0)[l]{$y$\strut}}
  \put(974.52,557.09){\fontsize{7.11}{8.54}\selectfont \makebox(75.0, 50.0)[l]{$X_2$\strut}}
  \put(1371.33,549.41){\fontsize{7.11}{8.54}\selectfont \makebox(25.0, 50.0)[l]{$X$\strut}}
  \put(1368.75,239.44){\fontsize{7.11}{8.54}\selectfont \makebox(25.0, 50.0)[l]{$Y$\strut}}
  \put(113.57,555.32){\fontsize{7.11}{8.54}\selectfont \makebox(75.0, 50.0)[l]{$X_0$\strut}}
  \put(392.28,562.47){\fontsize{7.11}{8.54}\selectfont \makebox(100.0, 50.0)[l]{$X_{11}$\strut}}
  \put(678.56,562.47){\fontsize{7.11}{8.54}\selectfont \makebox(100.0, 50.0)[l]{$X_{12}$\strut}}
  \put(395.46,243.06){\fontsize{7.11}{8.54}\selectfont \makebox(100.0, 50.0)[l]{$Y_{11}$\strut}}
  \put(686.81,239.44){\fontsize{7.11}{8.54}\selectfont \makebox(100.0, 50.0)[l]{$Y_{12}$\strut}}
  \put(118.75,250.60){\fontsize{7.11}{8.54}\selectfont \makebox(75.0, 50.0)[l]{$Y_0$\strut}}
\end{picture}%
\fi
\caption[Figure 1: ]{\mbox{ The bi-partition $(\widetilde{X},\widetilde{Y})$ of $V(G)$. }}\label{fig-1}
\end{figure}

By Claim \ref{claim-6} we have $e(B_{11},X_2)=0$. Using $e(B_0, X)=0$ and \eqref{eq_5.2}, we obtain that
\[
e(\widetilde{X}) = e(X_0)+e(X_{11})+e(X_{12})+e(X_0,X_{11}\cup X_{12}\cup X_2)+e(X_{11},X_{12})+e(X_{12},X_2)+e(X_{12},X).
\]
Similarly,
\[
e(\widetilde{Y}) = e(Y_0)+e(Y_{11})+e(Y_{12})+e(Y_0,Y_{11}\cup Y_{12}\cup Y_2)+e(Y_{11},Y_{12})+e(Y_{11},Y_2)+e(Y_{11},Y).
\]
Note that
\begin{align*}
&e(X_0)\leq \binom{|X_0|}{2},\ e(X_{11})\leq \binom{|X_{11}|}{2},\ e(X_{12})\leq \binom{|X_{12}|}{2},\\[5pt]
&e(X_0,X_{11}\cup X_{12}\cup X_2)\leq |X_0|(|X_{11}|+|X_{12}|+|X_2|),\\[5pt]
&e(X_{11},X_{12})\leq |X_{11}||X_{12}|,\ e(X_{12},X_2)\leq |X_{12}||X_2|.
\end{align*}
Similarly,
\begin{align*}
&e(Y_0)\leq \binom{|Y_0|}{2},\ e(Y_{11})\leq \binom{|Y_{11}|}{2},\ e(Y_{12})\leq \binom{|Y_{12}|}{2},\\[5pt]
&e(Y_0,Y_{11}\cup Y_{12}\cup Y_2)\leq |Y_0|(|Y_{11}|+|Y_{12}|+|Y_2|),\\[5pt]
 &e(Y_{11},Y_{12})\leq |Y_{11}||Y_{12}|,\ e(Y_{11},Y_2)\leq |Y_{11}||Y_2|.
\end{align*}
By the definitions of $B_{12}$ and $B_{11}$, $e(X_{12},X)\leq |X_{12}|$  and $e(Y_{11},Y)\leq |Y_{11}|$.

Recall that $f(x):= \binom{x}{2}-\left\lfloor\frac{x^2}{4}\right\rfloor+\left\lfloor\frac{x}{2}\right\rfloor$. Then
\begin{align*}
&e(X_0)+e(Y_0)\leq \binom{|X_0|}{2}+\binom{|Y_0|}{2}=\binom{z}{2} -\left\lfloor\frac{z^2}{4}\right\rfloor=f(z) -\left\lfloor\frac{z}{2}\right\rfloor,\\[5pt]
&e(X_{11})+e(Y_{11})+e(Y_{11},Y)\leq \binom{|X_{11}|}{2}+\binom{|Y_{11}|}{2}+|Y_{11}|=\binom{p}{2} -\left\lfloor\frac{p^2}{4}\right\rfloor+\left\lfloor\frac{p}{2}\right\rfloor=f(p),\\[3pt]
&e(X_{12})+e(Y_{12})+e(X_{12},X)\leq \binom{|X_{12}|}{2}+\binom{|Y_{12}|}{2}+|X_{12}|=\binom{q}{2} -\left\lfloor\frac{q^2}{4}\right\rfloor+\left\lfloor\frac{q}{2}\right\rfloor=f(q).
\end{align*}
Moreover,
\begin{align*}
&e(X_0,X_{11}\cup X_{12}\cup X_2)+e(X_{11},X_{12})\leq \left\lfloor\frac{z}{2}\right\rfloor\left(\left\lceil\frac{p}{2}\right \rceil+\left\lfloor\frac{q}{2}\right\rfloor+x\right)+\left\lceil\frac{p}{2}\right \rceil\left\lfloor\frac{q}{2}\right\rfloor,\\[3pt]
&e(Y_0,Y_{11}\cup Y_{12}\cup Y_2)+e(Y_{11},Y_{12})\leq \left\lceil\frac{z}{2}\right\rceil\left(\left\lfloor\frac{p}{2}\right\rfloor+\left\lceil\frac{q}{2}\right \rceil+y\right)+\left\lfloor\frac{p}{2}\right\rfloor\left\lceil\frac{q}{2}\right \rceil,
\end{align*}
and
\begin{align*}
&e(X_{12},X_2)\leq \left\lfloor\frac{q}{2}\right\rfloor x,\ e(Y_{11},Y_2)\leq \left\lfloor\frac{p}{2}\right\rfloor y.
\end{align*}
Thus,
\begin{align*}
 e(\widetilde{X})+e(\widetilde{Y}) &\leq f(p)+f(q)+f(z)- \left\lfloor\frac{z}{2}\right\rfloor+\left\lfloor\frac{z}{2}\right\rfloor\left(\left\lceil\frac{p}{2}\right \rceil+\left\lfloor\frac{q}{2}\right\rfloor+x\right)\nonumber\\[3pt]&+\left\lceil\frac{z}{2}\right\rceil\left(
 \left\lfloor\frac{p}{2}\right\rfloor+\left\lceil\frac{q}{2}\right \rceil+y\right)+\left\lceil\frac{p}{2}\right \rceil \left\lfloor\frac{q}{2}\right\rfloor+\left\lfloor\frac{p}{2}\right\rfloor\left\lceil\frac{q}{2}\right \rceil+\left\lfloor\frac{q}{2}\right\rfloor x+\left\lfloor\frac{p}{2}\right\rfloor y.
\end{align*}
Using \eqref{equalities2.1-1} and \eqref{equalities2.1-2}, we obtain that
\begin{align*}
 e(\widetilde{X})+e(\widetilde{Y}) \leq &f(p)+f(q)+f(z)- \left\lfloor\frac{z}{2}\right\rfloor+\left\lfloor\frac{zp}{2}\right\rfloor
 +\left\lceil\frac{zq}{2}\right \rceil+\left\lfloor\frac{z}{2}\right\rfloor x\\[3pt]
 & +\left\lceil\frac{z}{2}\right\rceil y+\left\lfloor\frac{pq}{2}\right\rfloor+\left\lfloor\frac{q}{2}\right\rfloor x+\left\lfloor\frac{p}{2}\right\rfloor y.
\end{align*}
By \eqref{equalities2.2} and \eqref{equalities2.4}, we arrive at
\begin{align*}
 e(\widetilde{X})+e(\widetilde{Y}) \leq &f(p)+f(q)+f(z)- \left\lfloor\frac{z}{2}\right\rfloor+\left\lfloor\frac{zp}{2}\right\rfloor
 +\left\lceil\frac{zq}{2}\right \rceil+\left\lfloor\frac{z x}{2}\right\rfloor\nonumber\\[3pt]
 & +\left\lceil\frac{zy}{2}\right\rceil+ \left\lfloor\frac{y}{2}\right\rfloor +\left\lfloor\frac{pq}{2}\right\rfloor+\left\lfloor\frac{q x}{2}\right\rfloor +\left\lfloor\frac{p y}{2}\right\rfloor.
\end{align*}
Since
\begin{align*}
f(p+q+z)-f(p)-f(q)-f(z)&\overset{\eqref{equalities2.5}}{=}f(p+q)-f(p)-f(q)
+\left\lceil\frac{z(p+q)}{2}\right\rceil\\[3pt]
&\overset{\eqref{equalities2.5}}{=} \left\lceil\frac{z(p+q)}{2}\right\rceil
+\left\lceil\frac{pq}{2}\right\rceil\\[3pt]
&\overset{\eqref{equalities2.3}}{\geq}\left\lfloor\frac{zp}{2}\right\rfloor
 +\left\lceil\frac{zq}{2}\right \rceil
+\left\lceil\frac{pq}{2}\right\rceil,
\end{align*}
we have
\begin{align*}
e(\widetilde{X})+e(\widetilde{Y}) &\leq f(p+q+z)- \left\lfloor\frac{z}{2}\right\rfloor
 +\left\lfloor\frac{z x}{2}\right\rfloor+\left\lceil\frac{zy}{2}\right\rceil+ \left\lfloor\frac{y}{2}\right\rfloor +\left\lfloor\frac{q x}{2}\right\rfloor +\left\lfloor\frac{p y}{2}\right\rfloor.
\end{align*}
Then,
\begin{align}
e(\widetilde{X})+e(\widetilde{Y})
&\leq f(p+q+z) +\left\lfloor\frac{zx}{2}\right\rfloor+\left\lceil\frac{zy}{2}\right\rceil+ \left\lfloor\frac{y}{2}\right\rfloor+\left\lfloor\frac{qx}{2}\right\rfloor+\left\lfloor\frac{py}{2}\right\rfloor
\label{ineq-3.3}\\[3pt]
&\overset{\eqref{equalities2.5}}{=} f(p+q+z+x)-f(x)-\left\lceil\frac{x(p+q+z)}{2}\right\rceil +\left\lfloor\frac{zx}{2}\right\rfloor+\left\lceil\frac{zy}{2}\right\rceil+ \left\lfloor\frac{y}{2}\right\rfloor\nonumber\\[3pt]
&\qquad +\left\lfloor\frac{qx}{2}\right\rfloor+\left\lfloor\frac{py}{2}\right\rfloor \nonumber\\[3pt]
&\overset{\eqref{equalities2.3}}{\leq}f(p+q+z+x)+\left\lceil\frac{zy}{2}\right\rceil+ \left\lfloor\frac{y}{2}\right\rfloor+\left\lfloor\frac{py}{2}\right\rfloor\label{ineq-3.5}\\[3pt]
&\overset{\eqref{equalities2.5}}{\leq} f(p+q+z+x+y)-f(y)-\left\lceil\frac{y(p+q+z+x)}{2}\right\rceil +\left\lceil\frac{zy}{2}\right\rceil+ \left\lfloor\frac{y}{2}\right\rfloor+\left\lfloor\frac{py}{2}\right\rfloor\nonumber\\[3pt]
&\leq  f(p+q+z+x+y)-\left(\left\lceil\frac{y(p+q+z+x)}{2}\right\rceil -\left\lceil\frac{zy}{2}\right\rceil-\left\lfloor\frac{py}{2}\right\rfloor\right)\label{ineq-3.6}\\[3pt]
&\overset{\eqref{equalities2.3}}{\leq} f(p+q+z+x+y)\nonumber.
\end{align}
Since $f(x)$ is an increasing function and $p+q+z+x+y \leq t+1$,
\[
\binom{\lfloor\frac{t+2}{2}\rfloor}{2}+\binom{\lceil\frac{t+2}{2}\rceil}{2}\leq \gamma_2(G)\leq e(\tilde{X})+e(\tilde{Y}) \leq f(t+1)\overset{\eqref{equalities2.6}}{=}\binom{\lfloor\frac{t+2}{2}\rfloor}{2}+\binom{\lceil\frac{t+2}{2}\rceil}{2}.
\]
Then we have equality in \eqref{ineq-3.3}, \eqref{ineq-3.5}, \eqref{ineq-3.6} and
\begin{align}\label{ineq-5.6}
p+q+z+x+y=t+1.
\end{align}
It follows that
\[
\left\lfloor\frac{z}{2}\right\rfloor=0,\
f(x)=0,\ f(y)=\left\lfloor\frac{y}{2}\right\rfloor.
\]
Hence $z\leq 1$, $x\leq 1$ and $y\leq 2$.

Recall that $(B_0, B_{11}, B_{12}, X_{2}, Y_{2})$ is partition of $B$ and $|B_0|=z\leq 1$. Let us define a different partition of $B_{11}$ and $B_{12}$. Let $X_{11}^*=\cup_{v\in Y_2} N(v,B_{11})$  and $Y_{12}^*=\cup_{u\in X_2} N(u,B_{12})$. Set $|B_{11}\setminus X_{11}^*|=p_1$, $|X_{11}^*|=p_2$, $|B_{12}\setminus Y_{12}^*|=q_1$ and $|Y_{12}^*|=q_2$. Clearly, $p=p_1+p_2$ and $q=q_1+q_2$.
 We  partition $B_{11}\setminus X_{11}^*$ into $X_{11}'$ and $Y_{11}'$ such that $|X_{11}'|=\left\lceil\frac{p_1}{2}\right\rceil$ and $|X_{12}'|=\left\lfloor\frac{p_1}{2}\right\rfloor$. Similarly, we  partition $B_{12}\setminus Y_{12}^*$ into $X_{12}'$ and $Y_{12}'$ such that $|X_{12}'|=\left\lfloor\frac{q_1}{2}\right\rfloor$ and $|Y_{12}'|=\left\lceil\frac{q_1}{2}\right\rceil$. Let
\[
X'= X_{11}'\cup X_{11}^*\cup X_{12}'\cup X_2\cup X\cup B_0,\  Y'=Y_{11}'\cup Y_{12}^*\cup Y_{12}'\cup Y_2\cup Y.
\]
Clearly, $(X',Y')$ is also a bi-partition of $V(G)$ (as shown in Figure \ref{fig-2}) and $\gamma_2(G)\leq e(X')+e(Y')$.
\begin{figure}[H]
\centering
\ifpdf
  \setlength{\unitlength}{0.05 mm}%
  \begin{picture}(1781.7, 903.1)(0,0)
  \put(0,0){\includegraphics{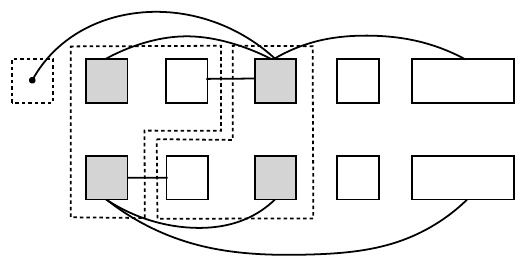}}
  \put(273.81,473.10){\fontsize{7.11}{8.54}\selectfont \makebox(125.0, 50.0)[l]{ $\left\lceil\frac{p_1}{2}\right \rceil $\strut}}
  \put(404.90,769.67){\fontsize{7.11}{8.54}\selectfont \makebox(100.0, 50.0)[l]{$B_{11}$\strut}}
  \put(856.20,99.18){\fontsize{7.11}{8.54}\selectfont \makebox(100.0, 50.0)[l]{$B_{12}$\strut}}
  \put(1170.72,267.46){\fontsize{7.11}{8.54}\selectfont \makebox(75.0, 50.0)[cc]{$Y_2$\strut}}
  \put(297.45,402.19){\fontsize{7.11}{8.54}\selectfont \makebox(125.0, 50.0)[l]{$\left\lfloor\frac{p_1}{2}\right\rfloor $\strut}}
  \put(602.61,486.00){\fontsize{7.11}{8.54}\selectfont \makebox(75.0, 50.0)[l]{$p_2 $\strut}}
  \put(600.46,384.99){\fontsize{7.11}{8.54}\selectfont \makebox(75.0, 50.0)[l]{$q_2$\strut}}
  \put(867.28,483.02){\fontsize{7.11}{8.54}\selectfont \makebox(125.0, 50.0)[l]{$\left\lfloor\frac{q_1}{2}\right\rfloor$\strut}}
  \put(869.80,398.60){\fontsize{7.11}{8.54}\selectfont \makebox(125.0, 50.0)[l]{$\left\lceil\frac{q_1}{2}\right \rceil$\strut}}
  \put(1183.52,477.98){\fontsize{7.11}{8.54}\selectfont \makebox(25.0, 50.0)[l]{$x$\strut}}
  \put(1191.08,393.56){\fontsize{7.11}{8.54}\selectfont \makebox(25.0, 50.0)[l]{$y$\strut}}
  \put(1169.96,590.25){\fontsize{7.11}{8.54}\selectfont \makebox(75.0, 50.0)[cc]{$X_2$\strut}}
  \put(1562.77,601.95){\fontsize{7.11}{8.54}\selectfont \makebox(25.0, 50.0)[l]{$X$\strut}}
  \put(1562.77,275.63){\fontsize{7.11}{8.54}\selectfont \makebox(25.0, 50.0)[l]{$Y$\strut}}
  \put(90.11,481.92){\fontsize{7.11}{8.54}\selectfont \makebox(25.0, 50.0)[l]{$z$\strut}}
  \put(87.88,723.66){\fontsize{7.11}{8.54}\selectfont \makebox(75.0, 50.0)[l]{$B_0$\strut}}
  \put(300.17,609.56){\fontsize{7.11}{8.54}\selectfont \makebox(125.0, 50.0)[cc]{$X_{11}'$\strut}}
  \put(309.63,279.41){\fontsize{7.11}{8.54}\selectfont \makebox(125.0, 50.0)[l]{$Y_{11}'$\strut}}
  \put(572.19,601.84){\fontsize{7.11}{8.54}\selectfont \makebox(125.0, 50.0)[cc]{$X_{11}^*$\strut}}
  \put(596.31,277.34){\fontsize{7.11}{8.54}\selectfont \makebox(125.0, 50.0)[l]{$Y_{12}^*$\strut}}
  \put(872.06,600.08){\fontsize{7.11}{8.54}\selectfont \makebox(125.0, 50.0)[cc]{$X_{12}'$\strut}}
  \put(876.09,273.76){\fontsize{7.11}{8.54}\selectfont \makebox(125.0, 50.0)[cc]{$Y_{12}'$\strut}}
  \end{picture}%
\else
  \setlength{\unitlength}{0.05 mm}%
  \begin{picture}(1781.7, 903.1)(0,0)
  \put(0,0){\includegraphics{5}}
  \put(273.81,473.10){\fontsize{7.11}{8.54}\selectfont \makebox(125.0, 50.0)[l]{ $\left\lceil\frac{p_1}{2}\right \rceil $\strut}}
  \put(404.90,769.67){\fontsize{7.11}{8.54}\selectfont \makebox(100.0, 50.0)[l]{$B_{11}$\strut}}
  \put(856.20,99.18){\fontsize{7.11}{8.54}\selectfont \makebox(100.0, 50.0)[l]{$B_{12}$\strut}}
  \put(1170.72,267.46){\fontsize{7.11}{8.54}\selectfont \makebox(75.0, 50.0)[cc]{$Y_2$\strut}}
  \put(297.45,402.19){\fontsize{7.11}{8.54}\selectfont \makebox(125.0, 50.0)[l]{$\left\lfloor\frac{p_1}{2}\right\rfloor $\strut}}
  \put(602.61,486.00){\fontsize{7.11}{8.54}\selectfont \makebox(75.0, 50.0)[l]{$p_2 $\strut}}
  \put(600.46,384.99){\fontsize{7.11}{8.54}\selectfont \makebox(75.0, 50.0)[l]{$q_2$\strut}}
  \put(867.28,483.02){\fontsize{7.11}{8.54}\selectfont \makebox(125.0, 50.0)[l]{$\left\lfloor\frac{q_1}{2}\right\rfloor$\strut}}
  \put(869.80,398.60){\fontsize{7.11}{8.54}\selectfont \makebox(125.0, 50.0)[l]{$\left\lceil\frac{q_1}{2}\right \rceil$\strut}}
  \put(1183.52,477.98){\fontsize{7.11}{8.54}\selectfont \makebox(25.0, 50.0)[l]{$x$\strut}}
  \put(1191.08,393.56){\fontsize{7.11}{8.54}\selectfont \makebox(25.0, 50.0)[l]{$y$\strut}}
  \put(1169.96,590.25){\fontsize{7.11}{8.54}\selectfont \makebox(75.0, 50.0)[cc]{$X_2$\strut}}
  \put(1562.77,601.95){\fontsize{7.11}{8.54}\selectfont \makebox(25.0, 50.0)[l]{$X$\strut}}
  \put(1562.77,275.63){\fontsize{7.11}{8.54}\selectfont \makebox(25.0, 50.0)[l]{$Y$\strut}}
  \put(90.11,481.92){\fontsize{7.11}{8.54}\selectfont \makebox(25.0, 50.0)[l]{$z$\strut}}
  \put(87.88,723.66){\fontsize{7.11}{8.54}\selectfont \makebox(75.0, 50.0)[l]{$B_0$\strut}}
  \put(300.17,609.56){\fontsize{7.11}{8.54}\selectfont \makebox(125.0, 50.0)[cc]{$X_{11}'$\strut}}
  \put(309.63,279.41){\fontsize{7.11}{8.54}\selectfont \makebox(125.0, 50.0)[l]{$Y_{11}'$\strut}}
  \put(572.19,601.84){\fontsize{7.11}{8.54}\selectfont \makebox(125.0, 50.0)[cc]{$X_{11}^*$\strut}}
  \put(596.31,277.34){\fontsize{7.11}{8.54}\selectfont \makebox(125.0, 50.0)[l]{$Y_{12}^*$\strut}}
  \put(872.06,600.08){\fontsize{7.11}{8.54}\selectfont \makebox(125.0, 50.0)[cc]{$X_{12}'$\strut}}
  \put(876.09,273.76){\fontsize{7.11}{8.54}\selectfont \makebox(125.0, 50.0)[cc]{$Y_{12}'$\strut}}
  \end{picture}%
\fi
\caption[Figure 2: ]{\mbox{ The bi-partition $(X',Y')$ of $V(G)$. }}\label{fig-2}
\end{figure}

\begin{claim}\label{claim-7}
If $X_{11}^*\neq \emptyset$, then $X_{11}^*$ is an independent set and $e(X_{11}', X_{11}^*)=0$. Similarly, if $Y_{12}^*\neq \emptyset$, then $Y_{12}^*$ is an independent set and $e(Y_{12}^*, Y_{12}')=0$.
\end{claim}
\begin{proof}
Suppose that there is an edge $wv$ with $w\in X_{11}^*$ and $v\in X_{11}'\cup X_{11}^*$. By the definition of $X_{11}^*$, there exists $u\in Y_2$ such that $uw\in E(G)$. Then $uwv$ be a path in $B$. Note that  $\deg(w,Y)= 1$, $\deg(v,Y)=1$ and $\deg(u,X)\geq 2$. It leads to a contradiction with Lemma \ref{lem-2.2} (vi). Thus the claim holds.
\end{proof}

\begin{claim}\label{claim-8}
 $p_2+q_2+x+y+z\leq 1$.
\end{claim}
\begin{proof}
If $z=1$, let $B_0=\{w\}$, then
by Lemma \ref{lem-2.2} (vi) we infer that  $e(\{w\},B_{11}\cup X_2)=0$ or $e(\{w\},B_{12}\cup Y_2)=0$. If $e(\{w\},B_{12}\cup Y_2)=0$ then one can interchange $X$ and $Y$, $X_2$ and $Y_2$, $X_{11}'$ and $Y_{12}'$, $X_{11}^*$ and $Y_{12}^*$,  $X_{12}'$ and $Y_{11}'$. Thus by symmetry we may assume that $e(\{w\},B_{11}\cup X_2)=0$. If $z=0$ then $e(B_0,B_{11}\cup X_2)=0$ and $e(B_0,B_{12}\cup Y_2)=0$ are obvious. Thus in either case we may assume that $e(B_0,B_{11}\cup X_2)=0$.

By Claim \ref{claim-7}, $e(X_{11}',X_{11}^*)=0$.  By Claim \ref{claim-6} we have $e(B_{11},X_2)=0$. By the definition of $Y_{12}^*$ we have $e(X_{12}',X_2)=0$. Using $e(B_0, X)=0$, \eqref{eq_5.2}  and by Claim \ref{claim-5}, Claim \ref{claim-7}, we arrive at
\[
e(X') =e(X_{11}')+e(X_{12}')+e(X_{11}',X_{12}')+e(X_{11}^*,X_{12}')+e(X_{12}',X)+e( B_0,X_{12}').
\]
Similarly,
$e(Y') = e(Y_{11}')+e(Y_{12}')+e(Y_{11}',Y_{12}^*)+e(Y_{11}',Y_{12}')+e(Y_{11}',Y)$.

Since $X_{11}',Y_{11}'$ is almost balanced partition of $p_1$ vertices, we infer that $e(X_{11}')+e(Y_{11}') \leq \binom{p_1}{2}-\lfloor\frac{p_1^2}{4}\rfloor$. Note that  each vertex in $Y_{11}'$ has exactly one neighbor in $Y$. Thus,
\begin{align}\label{two construction 1}
e(X_{11}')+e(Y_{11}')+e(Y_{11}',Y)\leq\binom{p_1}{2}-\left\lfloor\frac{p_1^2}{4}\right\rfloor +\left\lfloor\frac{p_1}{2}\right\rfloor = f(p_1).
\end{align}
Similarly,
\begin{align}\label{two construction 2}
e(X_{12}')+e(Y_{12}')+e(X_{12}',X)\leq f(q_1).
\end{align}
It is easy to see that
\begin{align}
&e(X_{11}',X_{12}')+e(Y_{11}',Y_{12}')\leq \left\lceil\frac{p_1}{2}\right \rceil \left\lfloor\frac{q_1}{2}\right\rfloor+\left\lfloor\frac{p_1}{2}\right\rfloor\left\lceil\frac{q_1}{2}\right\rceil
\overset{\eqref{equalities2.1-1}}{\leq }\left\lfloor\frac{p_1q_1}{2}\right\rfloor,\label{two construction 3}\\[3pt]
&e(X_{11}^*,X_{12}')\leq p_2\left\lfloor\frac{q_1}{2}\right\rfloor,\
e(Y_{11}',Y_{12}^*)\leq\left\lfloor\frac{p_1}{2}\right\rfloor q_2\label{two construction 4}
\end{align}
and $e(\{w\},X_{12}')\leq z\left\lfloor\frac{q_1}{2}\right\rfloor $. Adding \eqref{two construction 1}, \eqref{two construction 2},\eqref{two construction 3} and \eqref{two construction 4}, we get
\begin{align}
\gamma_2(G)&\leq e(X_{11}')+e(Y_{11}')+e(Y_{11}',Y)+e(X_{12}')+e(Y_{12}')+e(X_{12}',X)+e(X_{11}',X_{12}')\nonumber\\[3pt]
&+e(Y_{11}',Y_{12}')+e(X_{11}^*,X_{12}')+e(Y_{11}',Y_{12}^*)+e(\{w\},X_{12}')\nonumber\\[3pt]
&\leq f(p_1)+f(q_1)+\left\lfloor\frac{p_1q_1}{2}\right\rfloor
+ p_2\left\lfloor\frac{q_1}{2}\right\rfloor
+\left\lfloor\frac{p_1}{2}\right\rfloor q_2+z\left\lfloor\frac{q_1}{2}\right\rfloor\nonumber\\[3pt]
&\overset{\eqref{equalities2.5}}{\leq} f(p_1+q_1)+ p_2\left\lfloor\frac{q_1}{2}\right\rfloor
+\left\lfloor\frac{p_1}{2}\right\rfloor q_2+ z\left\lfloor\frac{q_1}{2}\right\rfloor\label{5.6}\\[3pt]
&\overset{\eqref{equalities2.5}}{\leq} f(p_1+q_1+p_2+q_2+x+y+z)-f(p_2+q_2+x+y+z)\nonumber\\[3pt]
&-\left\lceil\frac{(p_1+q_1)(p_2+q_2+x+y+z)}{2}\right\rceil
+ p_2\left\lfloor\frac{q_1}{2}\right\rfloor
+\left\lfloor\frac{p_1}{2}\right\rfloor q_2+ z\left\lfloor\frac{q_1}{2}\right\rfloor \nonumber
\end{align}
Since
\[
\left\lceil\frac{(p_1+q_1)(p_2+q_2+x+y+z)}{2}\right\rceil\overset{\eqref{equalities2.3}}{\geq} \left\lfloor\frac{q_1p_2}{2}\right\rfloor+\left\lfloor\frac{p_1 q_2}{2}\right\rfloor+\left\lfloor\frac{q_1z}{2}\right\rfloor\overset{\eqref{equalities2.2}}{\geq}  p_2\left\lfloor\frac{q_1}{2}\right\rfloor
+\left\lfloor\frac{p_1}{2}\right\rfloor q_2+ z\left\lfloor\frac{q_1}{2}\right\rfloor,
\]
we conclude that $\gamma_2(G)\leq f(t+1)-f(p_2+q_2+x+y+z)\leq f(t+1)$ and equality holds  if and only if $p_2+q_2+x+y+z\leq 1$.
\end{proof}

\begin{claim}\label{claim-9}
 $z=0$.
\end{claim}
\begin{proof}
By Claim \ref{claim-8}, $p_2+q_2+x+y+z\leq 1$. Suppose $z=1$. Then $p_2+q_2+x+y=0$. Since $G-B$ is $C_{2k+1}$-free, $e(G-B)\leq \left\lfloor\frac{(n-t-1)^2}{4}\right\rfloor$. Moreover, by \eqref{ineq-5.6} we infer $e(B,G-B)=p_1+q_1= t$ and $e(B)\leq \binom{t+1}{2}$. Thus
\[
e(G)\leq \left\lfloor\frac{(n-t-1)^2}{4}\right\rfloor+\binom{t+1}{2}+t<\left\lfloor\frac{(n-t-1)^2}{4}\right\rfloor+\binom{t+2}{2},
\]
contradicting \eqref{ineq-asmption1}.
\end{proof}
\begin{claim}\label{claim-10}
 $x+y=0$.
\end{claim}
\begin{proof}
 By Claim \ref{claim-8} we have $p_2+q_2+x+y+z\leq 1$. Suppose that $x+y=1$.  Then $p_2+q_2+z=0$. By \eqref{ineq-5.6}, $p_1+q_1= t$. Using \eqref{5.6}
we have
\begin{align*}
\gamma_2(G)\leq f(p_1+q_1)= f(t)< f(t+1),
\end{align*}
contradicting \eqref{ineq-asmption2}.
\end{proof}
By Claims \ref{claim-9} and \ref{claim-10},  $x=y=z=0$. Then $X_{11}^*=\cup_{v\in Y_2} N(v,B_{11})=\emptyset$  and $Y_{12}^*=\cup_{u\in X_2} N(u,B_{12})=\emptyset$, that is, $p_2=q_2=0$. Therefore $B=B_{11}\cup B_{12}$.

\begin{claim}\label{claim-11}
$G[B]$ is a clique.
\end{claim}
\begin{proof}
Since $e(B,X\cup Y)=e(B_{11},Y)+e(B_{12},X)\leq |B_{11}|+|B_{12}|=  t+1$. Moreover, $G[X\cup Y]$ is bipartite and $|X|+|Y| =n-t-1$, we have $e(X\cup Y) \leq \left\lfloor\frac{(n-t-1)^2}{4}\right\rfloor$. By $e(G)\geq \left\lfloor\frac{(n-t-1)^2}{4}\right\rfloor+\binom{t+2}{2}$,
we have
\[
e(B)\geq e(G)-e(B,X\cup Y)-e(X\cup Y)\geq  \binom{t+2}{2}-(t+1)=\binom{t+1}{2}.
\]
Thus $B$ forms a clique in $G$.
\end{proof}
\begin{claim}\label{claim-12}
Either $p=0$ or $q=0$ holds.
\end{claim}
\begin{proof}
Suppose that $p\geq 1$ and  $q \geq 1$. Then we show that  $|B|\geq 3$.  Otherwise, $|B_{11}|=|B_{12}|=1$ and it follows that
$G$ is a bipartite. That is, $\gamma_2(G)=0<f(t+1)$, contradicting \eqref{ineq-asmption2}. By Claim \ref{claim-11} and $p\geq 1$, $q\geq 1$, we can find a  path $uwv$ with  $u\in B_{11}$, $v\in B_{12}$ and $w\in B_{11}\cup B_{12}$. Since $\deg(u,Y)=1$, $\deg(v,X)=1$ and $\deg(w,X\cup Y)=1$, there exist $x\in X$ and $y\in Y$ such that $vx, uy\in E(G)$, contradicting Lemma \ref{lem-2.2} (vi). Thus the claim holds.
\end{proof}

Without loss generality, assume that $p=0$. By \eqref{ineq-5.6} and Claims \ref{claim-9}, \ref{claim-10} we obtain that  $q=t+1$. Similarly, if $q=0$, then $p=t+1$. That is, the equality \eqref{ineq-asmption1} holds if and only if $p=0$, $q=t+1$ or $p=t+1$, $q=0$. That is, $B=B_{11}$ or $B=B_{12}$. Recall that $B_{11}\cap B_{12}=\emptyset$. For each $uv\in E(B)$,  there exists a unique vertex $w_{uv}\in G-B$ such that $uw_{uv},vw_{uv}\in E(G)$. Otherwise, one can find distinct vertices $w_1, w_2$ in $X$ (or $Y$)  such that $uw_1, vw_2\in E(G)$, contradicting  Lemma \ref{lem-2.2} (v).  By Claim \ref{claim-11}, $B$ forms a clique in $G$. Thus, each vertex in $B$ is adjacent to the unique vertex in $G-B$. Hence $G$ is isomorphic to a subgraph of $H(n,t)$. By the assumption \eqref{ineq-asmption1} we infer that $G=H(n,t)$ up to isomorphism and the theorem is proven.
 \end{proof}

\section{Concluding remarks}

In this paper, we obtain an edge-stability result and a vertex-stability result for $C_{2k+1}$-free graphs with at least $\frac{(n-2k+1)^2}{4}+\binom{2k}{2}$ edges. It should be mention that Kor\'{a}ndi,  Roberts and  Scott \cite{KRS} proposed a  challenging conjecture concerning $C_{2k+1}$-free graphs with  $\frac{n^2}{4}-\delta n^2$ edges.

For a graph $F$, an $F$-blowup is a graph obtained  from $F$ by replacing each vertex with an independent set and replacing each edge by a complete bipartite graph. Let $\hb(n,k)$ be a  class of graphs consisting of all $C_{2k+3}$-blowups on $n$ vertices.

\begin{conj}[\cite{KRS}]\label{conj-1}
Fixed $k\geq 2$ and let $\delta>0$ be small enough. Then for any $\delta>\delta_0>0$ and large enough $n$, the following holds. For every $C_{2k+1}$-free graph $G$ on $n$ vertices with $(\frac{1}{4}-\delta_0)n^2\geq e(G) \geq (\frac{1}{4}-\delta)n^2$ edges, there is a graph $G^*\in \hb(n,k)$ satisfying $e(G^*)\geq e(G)$ and $\gamma_2(G^*)\geq \gamma_2(G)$.
\end{conj}

We define a class of graphs $\hg(n,k)$ on $n$ vertices as follows: a graph $G$ on $n$ vertices is in $\hg(n,k)$ if $G$ has exactly one block (2-connected component) being complete bipartite and all the other blocks are cliques of size at most $2k$.

Motivated by Theorem \ref{thm-main2} and Conjecture \ref{conj-1}, we make the following conjecture.

\begin{conj}\label{conj-2}
Fixed $k\geq 2$ and let $\delta>0$ be small enough. Then for any $\delta>0$ and large enough $n$, the following holds. For every $C_{2k+1}$-free graph $G$ on $n$ vertices with $e(G) \geq (\frac{1}{4}-\delta)n^2$ edges, there is a graph $G^*\in \hb(n,k)\cup \hg(n,k)$ satisfying $e(G^*)\geq e(G)$ and $\gamma_2(G^*)\geq \gamma_2(G)$.
\end{conj}

\end{document}